\documentclass{article}
\usepackage{amssymb, amsmath, amscd}

\newcommand \modk {{\cal M}_{\kappa}}

\newcommand \HH {{\cal H}}

\newcommand \la {\lambda}

\newcommand \Prob {\mathbb P}

\newcommand \oo {{\mathfrak o}}

\newtheorem {lemma} {Lemma}
\newtheorem {theorem} {Theorem}
\newtheorem {proposition} {Proposition}
\newtheorem{corollary}{Corollary}
\newtheorem{assumption} {Assumption}

\title{Finitely-additive measures on the asymptotic foliations of a Markov compactum.}
\author{Alexander I. Bufetov}
\date{} 
\begin{document}
\maketitle
\section{Introduction.}
\subsection{H{\"o}lder cocycles over translation flows.}
Let $\rho\geq 2$ be an integer, let $M$ be a compact orientable surface of genus $\rho$, 
and let $\omega$ be a holomorphic one-form on $M$. 
Denote by  
$
{\mathfrak m}=(\omega\wedge {\overline \omega})/2i
$
the area form induced by $\omega$ and assume that ${\mathfrak m}(M)=1$.

Let $h_t^+$ be the {\it vertical} flow on $M$ (i.e., the flow corresponding to $\Re(\omega)$); 
let $h_t^-$ be the {\it horizontal} flow on $M$ (i.e., the flow corresponding to $\Im(\omega)$).
The flows $h_t^+$, $h_t^-$ preserve the area ${\mathfrak m}$ and are uniquely ergodic. 

Take $x\in M$, $t_1, t_2\in {\mathbb R}_+$ and assume that the closure of the set
\begin{equation}
\label{admrectpsan}
\{h^+_{\tau_1} h^{-}_{\tau_2}x, 0\leq \tau_1< t_1, 0\leq \tau_2< t_2\}
\end{equation}
does not contain zeros of the form $\omega$.  Then the set (\ref{admrectpsan})
is called {\it an admissible rectangle} and denoted $\Pi(x, t_1, t_2)$.
Let ${\overline {\mathfrak C}}$ be the semi-ring of admissible rectangles.

Consider the linear space ${\cal Y}^+$ of H{\"o}lder cocyles $\Phi^+(x,t)$ over the vertical
flow $h_t^+$ which are invariant under horizontal holonomy. More precisely, a function 
$\Phi^+(x,t): M\times {\mathbb R}\to {\mathbb C}$ belongs to the space ${\cal Y}^+$
if it satisfies:
\begin{enumerate}
\item  $\Phi^+(x,t+s)=\Phi^+(x,t)+\Phi^+(h_t^+x, s)$;
\item  There exists $t_0>0$, $\theta>0$ such that
$|\Phi^+(x,t)|\leq t^{\theta}$ for all $x\in M$ and all $t\in {\mathbb R}$ satisfying $|t|<t_0$;
\item  If $\Pi(x, t_1, t_2)$ is an admissible rectangle, then
$\Phi^+(x, t_1)=\Phi^+(h_{t_2}^-x, t_1)$.
\end{enumerate}
For example, if a cocycle $\Phi_1^+$ is defined by $\Phi_1^+(x,t)=t$, then clearly $\Phi_1^+\in {\cal Y}^+$.

In the same way  define the space of ${\cal Y}^-$ of H{\"o}lder cocyles $\Phi^-(x,t)$ over 
the horizontal flow $h_t^-$  which are invariant under vertical holonomy, and set $\Phi_1^-(x,t)=t$.

Given $\Phi^+\in {\cal Y}^+$, $\Phi^-\in {\cal Y}^-$, a finitely additive
measure $\Phi^+\times \Phi^-$ on the semi-ring ${\overline {\mathfrak C}}$ of admissible rectangles
is introduced by the formula
\begin{equation}
\Phi^+\times \Phi^-(\Pi(x, t_1, t_2))=\Phi^+(x,t_1)\cdot \Phi^-(x, t_2).
\end{equation}

In particular, for $\Phi^-\in {\cal Y}^-$, set $m_{\Phi^-}=\Phi_1^+\times\Phi^-$:
\begin{equation}
\label{mphipsan}
m_{\Phi^-}(\Pi(x, t_1, t_2))=t_1\Phi^-(x, t_2).
\end{equation}
For any $\Phi^-\in {\cal Y}^-$ the measure 
$m_{\Phi^-}$ satisfies $(h_t^+)_*m_{\Phi^-}=m_{\Phi^-}$ and is an invariant distribution in the
sense of G.~Forni \cite{F1}, \cite{F2}. For instance, $m_{\Phi_1^-}={\mathfrak m}$.

A ${\mathbb C}$-linear pairing between ${\cal Y}^+$ and ${\cal Y^-}$ is given, for
$\Phi^+\in {\cal Y}^+$, $\Phi^-\in {\cal Y}^-$,
by the formula
\begin{equation}
\label{psanpairing}
< \Phi^+, \Phi^->=\Phi^+\times \Phi^-(M)
\end{equation}

The space of Lipschitz functions is not invariant under $h_t^+$, and a larger 
function space $Lip_w^+(M, \omega)$ of weakly Lipschitz functions is introduced as follows.
A bounded measurable function $f$ belongs to $Lip_w^+(M, \omega)$ if there exists a
constant $C$, depending only on $f$, such that for any admissible rectangle $\Pi(x, t_1, t_2)$ we
have
\begin{equation}
\label{weaklippsan}
\big| \int_0^{t_1} f(h_t^+x)dt -\int_0^{t_1}  f(h_t^+(h^-_{t_2}x)dt \big|\leq C.
\end{equation}  
Let $C_f$ be the infimum of all $C$ satisfying (\ref{weaklippsan}). We norm $Lip_w^+(X)$ by setting
$$
||f||_{Lip_w^+}=\sup_X f+C_f.
$$
By definition, the space  $Lip_w^+(M, \omega)$ contains all Lipschitz functions on $M$ and is
invariant under $h_t^+$. We denote by $Lip_{w,0}^+(M, \omega)$ the
subspace of  $Lip_w^+(M, \omega)$ of
functions whose integral with respect to ${\mathfrak m}$ is $0$.

\subsection{Flows along the stable foliation of a pseudo-Anosov diffeomorphism.}

Assume that  $\theta_1>0$ and a diffeomorphism $g:M\to M$ are such that
\begin{equation}
g^{*}(\Re(\omega))=\exp(\theta_1)\Re(\omega); \ g^{*}(\Im(\omega))=\exp(-\theta_1)\Im(\omega).
\end{equation}

The diffeomorphism $g$ induces a linear automorphism $g^*$ of the cohomology space $H^1(M, {\mathbb C})$. 
Denote by $E^+$ the expanding subspace of $g^*$
(in other words, $E^+$ is the subspace spanned by vectors 
corresponding to Jordan cells of $g^*$ with eigenvalues exceeding $1$ in absolute value).
The action of $g$ on ${\cal Y}^+$ is given by $g^*\Phi^+(x,t)=\Phi^+(gx, \exp(\theta_1)t).$                             

\begin{proposition}
\label{isomeplus}
There exists a $g^*$-equivariant isomorphism between $E^+$ and ${\cal Y}^+$.
\end{proposition}
\begin{theorem}
\label{multiplicpsan}
There exists a continuous mapping $\Xi^+: Lip_w^+(M, \omega)\to {\cal Y}^+$ such that 
for any $f\in Lip_w^+(M, \omega)$, any $x\in X$ and any $T>0$ we have
$$
\big| \int_0^T f\circ h_t^+(x) dt -\Xi^+(f)\big(x,T\big)\big|<C_{\varepsilon}||f||_{Lip_w^+}
(1+\log(1+T))^{2\rho+1}.
$$
The mapping $\Xi^+$ satisfies $\Xi^+(f\circ h^+_t)=\Xi^+(f)$ and $\Xi^+(f\circ g)=g^*\Xi^+(f)$. 
\end{theorem}
The mapping $\Xi^+$ is constructed as follows. 
By Proposition \ref{isomeplus} applied to the flow $h_t^-$, there exists a $g$-equivariant 
isomorphism between ${\cal Y}^-$ and the contracting space for the action of $g^*$ on  $H^1(M, 
{\mathbb C})$ (in other words, the subspace spanned by vectors
corresponding to Jordan cells with eigenvalues strictly less than $1$ in absolute value).
\begin{proposition}
The pairing $<, >$ given by (\ref{psanpairing}) is nondegenerate and $g^*$-invariant.
\end{proposition}
{\bf Remark.} Under the identification of ${\cal Y}^+$ and  ${\cal Y}^-$ 
with respective subspaces of $H^1(M, {\mathbb C})$, the pairing $<,>$ is taken to 
the cup-product on $H^1(M, {\mathbb C})$ (see Proposition 4.19 in Veech \cite{V3}).

If $f\in Lip_w^+(M, \omega)$, then $f$ is Riemann-integrable with 
respect to $m_{\Phi^-}$ for any $\Phi^-\in {\cal Y}^-$ (see (\ref{intfdpsi}) for 
a precise definition of the integral). 
Assign to $f$ a cocycle $\Phi_f^+$ in such a way that for all $\Phi^-\in {\cal Y}^-$ we have
\begin{equation}
\label{defphifpsan}
<\Phi_f^+, \Phi^->=\int_M f dm_{\Phi^-}.
\end{equation} 
By definition,  $\Phi^+_{f\circ h_t^+}=\Phi_f^+$. The mapping $\Xi^+$ of Theorem \ref{multiplicpsan} is given by the formula 
\begin{equation}
\label{xiphi}
\Xi^+(f)=\Phi^+_f.
\end{equation}

The first eigenvalue for the action of $g^*$ on $E^+$ is $\exp(\theta_1)$ and is always simple. 
If its second eigenvalue has the form $\exp(\theta_2)$, where $\theta_2>0$, and is simple as 
well, then the following limit theorem holds for $h_t^+$.

Given a bounded measurable function $f: X\to {\mathbb R}$ and $x\in X$,
introduce a continuous function ${\mathfrak S}_n[f,x]$ on the unit interval by the formula
\begin{equation}
\label{normsumspsan}
{\mathfrak S}_n[f,x](\tau)=\int_0^{\tau\exp(n\theta_1)} f\circ h^+_t(x)dt.
\end{equation}
The functions ${\mathfrak S}_n[f,x]$ are $C[0,1]$-valued random variables on the probability space $(M, {\mathfrak m})$.

\begin{theorem}
\label{limthmpsan}
If $g^*|_{E^+}$ has a simple, real second eigenvalue $\exp(\theta_2)$,
$\theta_2>0$, then there exists a continuous functional $\alpha: Lip_{w}^+(M, \omega)\to {\mathbb 
R}$ and a
compactly supported non-degenerate
measure $\eta$ on $C[0,1]$ such that for any $f\in Lip_{w,0}^+(M, \omega)$ satisfying 
$\alpha(f)\neq 0$
the sequence of random variables
$$
\frac{{\mathfrak S}_n[f,x]}{\alpha(f)\exp(n\theta_2)}
$$
converges in distribution to $\eta$ as $n\to\infty$.
\end{theorem}

The functional $\alpha$ is constructed explicitly as follows.
Under the assumptions of the theorem the action of $g^*$ on $E^-$ has a simple eigenvalue $\exp(-\theta_2)$; 
let $v(2)$ be the eigenvector with eigenvalue $\exp(-\theta_2)$, let $\Phi_2^-\in {\cal Y}^-$ 
correspond to $v(2)$ by Proposition \ref{isomeplus} and $m_{\Phi_2^-}$ be given by (\ref{mphipsan}); 
then 
$$
\alpha(f)=\int fdm_{\Phi_2^-}.
$$

\subsection{Generic translation flows.}

Let $\rho\geq 2$ and let $\kappa=(\kappa_1, \dots, \kappa_{\sigma})$ be
a nonnegative integer vector such that  $\kappa_1+\dots+\kappa_{\sigma}=2\rho-2$.
Denote by $\modk$ the moduli space
of Riemann surfaces of genus $\rho$ endowed with a holomorphic differential
of area $1$ with singularities of orders $k_1, \dots, k_{\sigma}$
(the {\it stratum} in the moduli space of holomorphic differentials), and
let $\HH$ be a connected component of $\modk$.
Denote by $g_t$ the Teichm{\"u}ller flow on $\HH$
(see \cite{F2}, \cite{kz1}), and let ${\mathbb A}(t,X)$ be the Kontsevich-Zorich cocycle over $g_t$ \cite{kz1}.

Let $\Prob$ be a $g_t$-invariant ergodic probability
measure on $\HH$. For $X\in\HH$, $X=(M, \omega)$, let ${\cal Y}_X^+$, ${\cal Y}_X^-$ be the corresponding spaces of H{\"o}lder 
cocycles. Denote by $E_X^+$
the space spanned by the positive Lyapunov exponents of the Kontsevich-Zorich cocycle.

\begin{proposition}
For $\Prob$-almost all $X\in\HH$, we have 
$\dim {\cal Y}_X^+=\dim {\cal Y}_X^-=\dim E_X^+$, and 
the pairing $<,>$ between  ${\cal Y}_X^+$ and  ${\cal Y}_X^-$ is non-degenerate.
\end{proposition}
{\bf Remark.} In particular, if $\Prob$ is the Masur-Veech ``smooth" measure \cite{masur, V1}, then 
$\dim {\cal Y}_X^+=\dim {\cal Y}_X^-=\rho$.

Assign to $f\in Lip_w^+(M,\omega)$ a cocycle $\Phi_f^+$ by (\ref{defphifpsan}).
\begin{theorem}
\label{multiplicmoduli}
For any $\varepsilon>0$ there exists a constant $C_{\varepsilon}$ depending only on $\Prob$ such that for
$\Prob$-almost every $X\in\HH$, any $f\in Lip_w^+(X)$, any $x\in X$ and any $T>0$ we have
$$
\big| \int_0^T f\circ h_t^+(x) dt -\Phi^+_f(x,T)\big|<C_{\varepsilon}||f||_{Lip_w^+}(1+T^{\varepsilon}).
$$
\end{theorem}

If both the first and the second Lyapunov exponent of the measure $\Prob$ are positive and simple
(as, by the Avila-Viana Theorem \cite{AV},
is the case with the Masur-Veech ``smooth" measure on $\HH$), then the following limit theorem holds.

As before, consider a $C[0,1]$-valued random variable ${\mathfrak S}_t[f,x]$ on
$(M, {\mathfrak m})$ defined by the formula
$$
{\mathfrak S}_s[f,x](\tau)=\int_0^{\tau\exp(s)} f\circ h^{+}_t(x)dt.
$$
Let $||v||$ be the Hodge norm in $H^1(M, {\mathbb R})$.
Let $\theta_2>0$ be the second Lyapunov exponent of  the Kontsevich-Zorich cocycle
and let
$v_2(X)$ be a Lyapunov vector corresponding to $\theta_2$ (by our assumption, such a vector
is unique up to scalar multiplication). Introduce a real-valued multiplicative cocycle $H_2(t, X)$
over $g_t$ by the formula
\begin{equation}
H_2(t, X)=\frac{||A(t,X)v_2(X)||}{||v_2(X)||}.
\end{equation}

\begin{theorem}
\label{markovlimitmoduli}
Assume that
both the first and the second Lyapunov exponent of  the Kontsevich-Zorich cocycle with respect
to the measure $\Prob$ are positive and simple.
Then for $\Prob$-almost any $X^{\prime}\in\HH$ there exists
a non-degenerate compactly supported  measure $\eta_{X^{\prime}}$ on $C[0,1]$ and,
for $\Prob$-almost all $X, X^{\prime}\in\HH$, there
exists a sequence of moments
$s_n=s_n(X, X^{\prime})$ such that the following holds.
For $\Prob$-almost every $X\in\HH$ there exists a continuous functional
$$
{\mathfrak a}^{(X)}: Lip_w^+(X)\to {\mathbb R}
$$
such that for $\Prob$-almost every $X^{\prime}$ and
for any real-valued $f\in Lip_{w,0}^+(X)$
satisfying ${\mathfrak a}^{(X)}(f) \neq 0$, the sequence of $C[0,1]$-valued
random variables
$$
\frac{{\mathfrak S}_{s_n}[f,x](\tau)}{\big({\mathfrak a}^{(X)}(f)\big)H_2(s_n, X)}
$$
converges in distribution to $\eta_{X^{\prime}}$ as $n\to\infty$.
\end{theorem}

{\bf Acknowledgements.}
W. A.~Veech made the suggestion that G. Forni's invariant distributions for the vertical flow 
should admit a description in terms of cocycles for the horizontal flow, and I am greatly 
indebted to him. The observation that cocycles are dual objects to invariant distributions was made by 
G.~Forni, and I am deeply grateful to him. I am deeply grateful to H. Nakada who pointed out the reference
to S.~Ito's work \cite{ito} to me. I am deeply grateful to J.~Chaika, P.~Hubert, Yu.S.~Ilyashenko, 
E.~Lanneau and R.~Ryham for many helpful suggestions on improving the presentation.
I am deeply grateful to A.~Avila, X.~Bressaud, B.M.~Gurevich, A.V.~Klimenko, V.I.~Oseledets,
Ya.G.~Sinai, I.V.~Vyugin, J.-C.~Yoccoz for useful discussions. During the work on this paper, 
I was supported in part by the National
Science Foundation under grant DMS~0604386 and by the Edgar Odell Lovett Fund at Rice University.

\section{Asymptotic foliations of a Markov compactum.}
\subsection{Definitions and notation.}
Let $m\in {\mathbb N}$ and let $\Gamma$ be an oriented graph with $m$ vertices $\{1,\dots,m\}$
and possibly multiple
edges. We assume that that for each vertex there is an edge starting
from it and an edge ending in it.

Let ${\cal E}(\Gamma)$ be the set of edges of $\Gamma$.
For $e\in {\cal E}(\Gamma)$ we denote by $I(e)$ its initial vertex and by $F(e)$ its terminal
vertex. Let  $Q$ be the incidence matrix of $\Gamma$ defined by the formula
$$
Q_{ij}=\# \{e\in {\cal E}(\Gamma): I(e)=i, F(e)=j\}.
$$
By assumption, all entries of the matrix $Q$ are positive.
A finite word $e_1 \dots e_k$, $e_i\in  {\cal E}(\Gamma)$, will be called 
{\it admissible} if  $F(e_{i+1})=I(e_i)$, $i=1, \dots, k$. 

To the graph $\Gamma$ we assign a {\it Markov compactum} $X_{\Gamma}$, the space 
of bi-infinite paths along the edges: 
$$
X_{\Gamma}=\{x=\dots x_{-n}\dots x _0\dots x_n \dots, \  x_n\in {\cal E}(\Gamma), F(x_{n+1})=I(x_n) \}.
$$

{\bf Remark.} As $\Gamma$ will be fixed throughout this section, we shall often omit the subscript 
$\Gamma$ from notation 
and only insert it when the dependence on $\Gamma$ is underlined.

Cylinders in $X_{\Gamma}$ are subsets of the form $\{x: x_{n+1}=e_1, \dots, x_{n+k}=e_{k}\}$, where 
$n\in {\mathbb Z}$, $k\in {\mathbb N}$ and $e_1\dots e_k$ is an admissible word.
The family of all cylinders forms a semi-ring which we denote by ${\mathfrak C}$. 

For $x\in X$, $n\in {\mathbb Z}$, introduce the sets
$$
\gamma^+_n(x)=\{x^{\prime}\in X_{\Gamma}: x^{\prime}_t=x_t, t\geq n\};\ 
\gamma^-_n(x)=\{x^{\prime}\in X_{\Gamma}: x^{\prime}_t=x_t, t\leq n\};
$$
$$
\gamma^+_{\infty}(x)=\bigcup_{n\in {\mathbb Z}} \gamma_n^+(x);\ 
\gamma^-_{\infty}(x)=\bigcup_{n\in {\mathbb Z}} \gamma_n^-(x).
$$
The sets $\gamma^+_{\infty}(x)$ are leaves of the asymptotic foliation  ${\cal F}^+$ on the space $X_{\Gamma}$; 
the sets $\gamma^+_{\infty}(x)$ are leaves of the asymptotic foliation  ${\cal F}^-$ on $X_{\Gamma}$.

For $n\in {\mathbb Z}$ let ${\mathfrak C}^+_n$ be the collection of all subsets of $X_{\Gamma}$ of the form 
$\gamma^+_n(x)$, $n\in {\mathbb Z}$, $x\in X$; similarly, 
${\mathfrak C}^-_n$ is the collection of all subsets of the form  $\gamma^-_n(x)$. 
Set 
\begin{equation}
\label{cplus}
{\mathfrak C}^+=\bigcup\limits_{n\in {\mathbb Z}} {\mathfrak C}^+_n; \ 
{\mathfrak C}^-=\bigcup\limits_{n\in {\mathbb Z}} {\mathfrak C}^-_n. 
\end{equation}

The collection ${\mathfrak C}^+_n$ is a semi-ring for any $n\in 
{\mathbb Z}$.
Since every element of ${\mathfrak C}^+_n$ is a disjoint union of elements of ${\mathfrak C}^+_{n+1}$, 
the collection ${\mathfrak C}^+$
is a semi-ring as well. The same  statements hold for ${\mathfrak C}^-_n$ and  ${\mathfrak C}^-$.   

Let $\exp({\theta_1})$ be the spectral radius of the matrix $Q$, and let 
$h=(h_1, \dots, h_m)$ be the unique positive eigenvector of $Q$: we thus have $Qh=\exp(\theta_1)h$. 
Let $\la=(\la_1, \dots, \la_m)$ be the positive eigenvector of the transpose matrix $Q^t$: we have $Q^t\la=\exp(\theta_1)\la$.
The vectors $\la, h$ are normalized as follows:
\begin{equation}
\label{lahnorm}
\sum_{i=1}^m \la_i=1; \ \sum_{i=1}^m \la_ih_i=1.
\end{equation}

Introduce a sigma-additive positive measure $\Phi_1^+$ on the semi-ring 
${\mathfrak C}^+$ by the formula 
\begin{equation}
\label{defPhioneplus}
\Phi_1^{+}(\gamma_n^+(x))=h_{F(x_n)}\exp((n-1)\theta_1)
\end{equation}
and a sigma-additive positive measure $\Phi_1^-$ on the semi-ring
${\mathfrak C}^-$ by the formula
\begin{equation}
\label{defPhioneminus}
\Phi_1^-(\gamma_n^-(x))=\la_{I(x_n)}\exp(-n\theta_1).
\end{equation}

Let  $n\in {\mathbb Z}$, $k\in{\mathbb N}$, and let $e_1\dots e_k$ be an admissible word.
The Parry measure $\nu$ on $X_{\Gamma}$ is defined by the formula
\begin{equation}
\label{defnu}
\nu(\{x: x_{n+1}=e_1, \dots, x_{n+k}=e_{k}\})=\la_{I(e_k)}h_{F(e_1)}\exp(-k\theta_1).
\end{equation}

The measures $\Phi_1^+$, $\Phi_1^-$ are conditional measures of the Parry measure $\nu$ in the following sense.
If $C\in{\mathfrak C}$, then $\gamma_{\infty}^+(x)\cap C\in {\mathfrak C}^+$,  $\gamma_{\infty}^-(x)\cap C\in {\mathfrak C}^-$ 
for any $x\in C$, and we have
\begin{equation}
\label{condmeas}
\nu(C)=\Phi_1^+(\gamma_{\infty}^+(x)\cap C)\cdot \Phi_1^-(\gamma_{\infty}^-(x)\cap C).
\end{equation}

\subsection{Finitely-additive measures on leaves of asymptotic foliations.}

Given $v \in {\mathbb C}^m$, write  
\begin{equation}
\label{lonenorm}
|v|=\sum\limits_{i=1}^m |v_i|.
\end{equation}
The norms of all matrices in this paper are understood with respect to this norm.  
Consider the direct-sum decomposition
$$
{\mathbb C}^m=E^+\oplus E^-, 
$$
where $E^+$ is spanned
by Jordan cells of eigenvalues of $Q$ with absolute value exceeding $1$, and $E^-$ is spanned by Jordan
cells corresponding to eigenvalues of $Q$ with absolute value at most $1$. 
Let $v\in E^+$ and for all $n\in {\mathbb Z}$ set $v^{(n)}=Q^{n}v$
(note that $Q|_{E^+}$ is by definition invertible).
Introduce a
finitely-additive complex-valued measure $\Phi^+_v$ on the semi-ring ${\mathfrak C}^+$ 
(defined in (\ref{cplus}))
by the formula
\begin{equation}
\label{defPhiplusv}
\Phi^{+}_v(\gamma_{n+1}^+(x))=(v^{(n)})_{F(x_{n+1})}.
\end{equation}

The measure  $\Phi^+_v$ is invariant under holonomy along ${\cal F}^-$: by definition, we have the following  
\begin{proposition}
\label{holoplus}
If $F(x_n)=F(x^{\prime}_n)$, then $\Phi^{+}_v(\gamma_n^+(x))=\Phi^{+}_v(\gamma_n^+(x^{\prime}))$.
\end{proposition}

The measures $\Phi^+_v$ span a complex linear space, which we denote  ${\cal Y}^+$ (or, sometimes,  ${\cal Y}^+_{\Gamma}$, when dependence on 
$\Gamma$ is stressed.) The map 
\begin{equation} 
\label{cali}
{\cal I}: v \to \Phi^+_v
\end{equation}
is an isomorphism between $E^+$ and  ${\cal Y}^+_{\Gamma}$.

For $Q^t$, we have the direct-sum decomposition
$$
{\mathbb C}^m={\tilde E}^+\oplus {\tilde E}^-,
$$
where ${\tilde E}^+$ is spanned
by Jordan cells of eigenvalues of $Q^t$ with absolute value exceeding $1$, and ${\tilde E}^-$ is spanned by
Jordan
cells corresponding to eigenvalues of $Q^t$ with absolute value at most $1$.
As before, for ${\tilde v}\in {\tilde E}^+$  set ${\tilde v}^{(n)}=(Q^t)^{n}{\tilde v}$
for all $n\in {\mathbb Z}$, and  introduce a
finitely-additive complex-valued measure $\Phi^-_{{\tilde v}}$ on the semi-ring ${\mathfrak C}^-$ 
(defined in (\ref{cplus})) by the formula
\begin{equation}
\label{defPhiminusv}
\Phi^{-}_{{\tilde v}}(\gamma_n^-(x))=({\tilde v}^{(-n)})_{I(x_n)}.
\end{equation}

By definition, the measure  $\Phi^-_{\tilde v}$ is invariant under holonomy along ${\cal F}^+$:
more precisely, we have the following
\begin{proposition}
\label{holominus}  
If $I(x_n)=I(x^{\prime}_n)$, then $\Phi^{-}_{{\tilde v}}(\gamma_n^-(x))=\Phi^{-}_{{\tilde v}}(\gamma_n^-(x^{\prime}))$.  
\end{proposition}

Let  ${\cal Y}^-_{\Gamma}$ be the space spanned by the measures
$\Phi^-_v$, $v\in {\tilde E}^+$.
The map 
\begin{equation}
\label{tildecali}
{\tilde {\cal I}}: v \to \Phi^-_v
\end{equation}
 is an isomorphism
between ${\tilde E}^+$ and ${\cal Y}^-_{\Gamma}$.

Let $\sigma: X_{\Gamma}\to X_{\Gamma}$ be the shift defined by $(\sigma x)_i=x_{i+1}$.
The shift $\sigma$ naturally acts on the spaces ${\cal Y}^+_{\Gamma}$, ${\cal Y}^-_{\Gamma}$:
given $\Phi\in {\cal Y}^+_{\Gamma}$ (or ${\cal Y}^-_{\Gamma}$), the measure $\sigma_*\Phi$ is 
defined, for $\gamma\in {\mathfrak C}^+$, by the formula
$$
\sigma_*\Phi(\gamma)=\Phi(\sigma\gamma).
$$
From the definitions we obtain
\begin{proposition}
\label{propcali}
The following  diagrams are commutative: 
$$
\begin{CD}
\label{shifteplus}
E^+@ >{\cal I}>> {\cal Y}^+_{\Gamma} \\
@VVQV           @AA{\sigma^*}A   \\
E^+@ >{\cal I}>> {\cal Y}^+_{\Gamma} \\
\end{CD}
$$
$$
\begin{CD}
\label{shiftetilde}
{\tilde E}^+@ >{\tilde {\cal I}}>> {\cal Y}^-_{\Gamma} \\
@VVQ^tV           @VV{\sigma^*}V   \\
{\tilde E}^+@ > {\tilde {\cal I}}>> {\cal Y}^-_{\Gamma} \\
\end{CD}
$$
\end{proposition}
\subsection{Pairings.}
\label{secpair}

Given $\Phi^+\in {\cal Y}^+$, $\Phi^-\in {\cal Y}^-$, introduce, in analogy with (\ref{condmeas}),  
a finitely additive measure $\Phi^+\times \Phi^-$ on the semi-ring ${\mathfrak C}$ 
of cylinders in $X_{\Gamma}$: for any $C\in {\mathfrak C}$ and $x\in  C$, set 
\begin{equation}
\label{prodmeas}
\Phi^+\times \Phi^-(C)=\Phi^+(\gamma^+_{\infty}(x)\cap C)\cdot \Phi^-(\gamma^-_{\infty}(x)\cap  C).
\end{equation}
Note that by Propositions \ref{holoplus}, \ref{holominus},
the right-hand side in (\ref{prodmeas}) does not depend on $x\in { C}$.

More explicitly, let  $v\in E^+$, ${\tilde v}\in {\tilde E}^+$, $\Phi^+_v={\cal I}(v)$,
$\Phi^-_{\tilde v}={\tilde {\cal I}}({\tilde v})$. As above, denote  $v^{(n)}=Q^{n}v$, ${\tilde v}^{(n)}=(Q^t)^{n}v$.
Let $n\in {\mathbb Z}$, $k\in {\mathbb N}$ 
and let $e_1\dots e_k$ be an admissible word. Then
\begin{equation}
\label{explicitproduct}
\Phi^+_v\times \Phi^-_{\tilde v}(\{x: x_{n+1}=e_1, \dots, x_{n+k}=e_{k}\})=
\big(v^{(n)}\big)_{F(e_1)}\big({\tilde v}^{(-n-k)}\big)_{I(e_{n+k})}.
\end{equation}

There is a natural ${\mathbb C}$-linear pairing $<,>$ between the spaces ${\cal Y}^+_{\Gamma}$ and
${\cal Y}^-_{\Gamma}$: for $\Phi^+\in {\cal Y}^+_{\Gamma}$, $\Phi^-\in {\cal Y}^-_{\Gamma}$, set
\begin{equation}
<\Phi^+, \Phi^->=\Phi^+\times \Phi^-(X_{\Gamma}).
\end{equation}
From (\ref{explicitproduct}) we derive
\begin{proposition}
Let  $v\in E^+$, ${\tilde v}\in {\tilde E}^+$, $\Phi^+_v={\cal I}_{\Gamma}(v)$,
$\Phi^-_{\tilde v}={\tilde {\cal I}}_{\Gamma}({\tilde v})$. Then
\begin{equation}
<\Phi^+_v, \Phi^-_{{\tilde v}}>=\sum_{i=1}^m v_i{\tilde v}_i.
\end{equation}
In particular, the pairing $<,>$ is non-degenerate and $\sigma^*$-invariant.
\end{proposition}
In particular, for $\Phi^-\in {\cal Y}^-$ denote 
\begin{equation}
\label{mphi}
m_{\Phi^-}=\Phi_1^+\times \Phi^-.
\end{equation}


\subsection{Weakly Lipschitz Functions.}
\label{secweaklip}

Introduce a function space  $Lip_w^+(X)$ in the following way.
A bounded Borel-measurable function $f:X \to {\mathbb C}$ belongs to the space $Lip_w^+(X)$ if there exists a
constant $C>0$ such that for all $n\geq 0$ and any $x, x^{\prime} \in X$ satisfying
$F(x_{n+1})=F(x^{\prime}_{n+1})$, we have
\begin{equation}
\label{weaklip} 
|\int_{\gamma_n^+(x)}fd\Phi_1^+  - \int_{\gamma_n^+(x^{\prime})} fd\Phi_1^+|\leq C.
\end{equation}
If $C_f$ be the infimum of all $C$ satisfying (\ref{weaklip}), then we norm $Lip_w^+(X)$ by setting
$$
||f||_{Lip_w^+}=\sup_X f+C_f.
$$
As before, let $Lip_{w,0}^+(X)$ be the subspace of $Lip_w^+(X)$ of functions 
whose integral with respect to $\nu$ is zero.

Take $\Phi^-\in {\cal Y}^-$. Any function $f\in Lip_w^+(X)$ is integrable with respect
to the measure $m_{\Phi^-}$, defined by (\ref{mphi}), in the following sense. Let
${\tilde v}\in E^-$ be the vector corresponding to $\Phi^-$ by (\ref{defPhiminusv}) and let
${\tilde v}^{(n)}=(Q^t)^n{\tilde v}$. Recall that
\begin{equation}
\label{tildevexpfast}
|{\tilde v}^{(-n)}|\to 0 \ {\rm exponentially \ fast \ as} \  n\to\infty.
\end{equation}

Take arbitrary points $x^{(n)}_i\in X$,  $n\in {\mathbb N}$ satisfying
\begin{equation}
\label{xin}
F((x_i^{(n)})_n)=i, \ i=1, \dots, m.
\end{equation}
and consider the expression
\begin{equation}
\label{intfdpsi}
\sum\limits_{i=1}^m \big( \int_{\gamma_n^+(x_i^{(n)})} fd\Phi_1^+  \big)\cdot
\big({\tilde v}^{(1-n)}\big)_i.
\end{equation}

By (\ref{weaklip}) and (\ref{tildevexpfast}), as $n\to\infty$ the expression (\ref{intfdpsi})
tends to a limit which does not depend on the particular choice of  $x^{(n)}_i$ satisfying (\ref{xin}).
This limit is denoted
$$
m_{\Phi^-}(f)=\int_X fdm_{\Phi^-}.
$$

Introduce a measure $\Phi_f^+\in {\cal Y}^+$ by requiring that for any
$\Phi^-\in {\cal Y}^-$ we have   
\begin{equation}
\label{phifplus}
<\Phi_f^+, \Phi^->=\int_X fdm_{\Phi^-}.
\end{equation}

Note that the mapping ${\Xi}^+:  Lip_w^+(X)\to {\cal Y}^+$ given by 
${\Xi}^+(f)=\Phi^+_f$ is continuous by definition and satisfies
\begin{equation}
\label{invarphiplus}
{\Xi}^+(f\circ \sigma)=\sigma^*{\Xi}^+(f).
\end{equation}

From the definitions we also have 
\begin{proposition}
Let $\Phi^+(1), \dots, \Phi^+(r)$ be a basis in ${\cal Y}^+$ 
and let $\Phi^-(1), \dots, \Phi^-(r)$ be the dual basis in ${\cal Y}^-$ with respect to the pairing $<,>$.
Then  for any $f\in Lip_w^+(X)$ we have 
$$
\Phi_f^+=\sum_{i=1}^r \big(m_{\Phi^-(i)}(f)\big)\Phi^+(i). 
$$
\end{proposition}

\subsection{Approximation.}

Let $\Theta$ be a finitely-additive complex-valued measure on the semi-ring ${\mathfrak C}^+_0$.
Assume that there exists a constant $\delta({\Theta})$  such that for all $x, x^{\prime} \in X$ and all
$n\geq 0$ we have
\begin{equation} 
|\Theta(\gamma_n^+(x))-\Theta(\gamma_n^+(x^{\prime}))|\leq \delta({\Theta}) \ {\rm if }
\ F(x_{n+1})=F(x^{\prime}_{n+1}).
\end{equation}

In this case $\Theta$ will be called a {\it weakly Lipschitz measure}.

\begin{lemma}
\label{thetapprox}
There exists a constant  $C_{\Gamma}$ depending only on $\Gamma$ such that the following is true.
Let $\Theta$ be a weakly Lipschitz  finitely-additive complex-valued measure on the semi-ring ${\mathfrak C}^+_0$.
Then there exists a unique $\Phi^+\in {\cal Y}^+_{\Gamma}$ such that for all $x\in X$ and
all $n>0$ we have
\begin{equation}
|\Theta(\gamma_n^+(x))-\Phi^+(\gamma_n^+(x))|\leq C_{\Gamma} \delta(\Theta)n^{m+1}.
\end{equation}
\end{lemma}

Assign to the graph $\Gamma$ the Markov compactum $Y_{\Gamma}$ of one-sided infinite
sequences of edges:
$$
Y=\{y=y_1\dots y_n \dots : y_n\in {\cal E}(\Gamma), F(y_{n+1})=I(y_n) \},
$$
and, as before, let $\sigma$ be the shift on $Y_{\Gamma}$: $(\sigma y)_i=y_{i+1}$.
For $y, y^{\prime}\in Y_{\Gamma}$, write $y^{\prime}\searrow y$ if $\sigma y^{\prime}=y$.

Lemma \ref{thetapprox} will be derived from
\begin{lemma}
\label{phinapprox}
There exists a constant  $C_{\Gamma}$ depending only on $\Gamma$ such that the following is true.
Let $\varphi_n$ be a sequence of  measurable complex-valued functions on  $Y_{\Gamma}$.
Assume that there exists a constant $\delta$ such that for all $y\in Y$ and all
$n\geq 0$ we have
\begin{equation}
\label{ysearrow}
|\varphi_{n+1}(y))-\sum\limits_{y^{\prime}\searrow y} \varphi_n(y^{\prime})|\leq \delta
\end{equation}
and for all $n\geq 0$ and all $y, {\tilde y}\in Y_{\Gamma}$ satisfying $F(y_1)=F({\tilde y}_1)$ we
have
\begin{equation}
\label{ytildey}
|\varphi_{n}(y))-\varphi_n({\tilde y})|\leq  \delta.
\end{equation}
Then there exists a unique $v\in E^+$ such that for all $y\in Y$ and all
$n>0$ we have
\begin{equation}
|\varphi_{n}(y))-(Q^n v)_{F(y_{n+1})}|\leq C_{\Gamma} \delta n^{m+1}.
\end{equation}
\end{lemma}
Proof of Lemma \ref{phinapprox}. Take arbitrary points $y(i)\in Y_{\Gamma}$ in such a way
that $$
F(y(i)_1)=i.
$$ 
Introduce a sequence of vectors $v(n)\in {\mathbb C}^m$ by the formula
$$
v(n)_i=\varphi_n(y(i)).
$$
From (\ref{ytildey}) for any $y\in Y$ we have
$$
|\varphi_n(y)-v(n)_{F(y_1)}|\leq \delta,
$$
and from (\ref{ysearrow}), (\ref{ytildey}) we have
$$
|Qv(n)-v(n+1)|\leq  \delta\cdot ||Q||.
$$

To prove Lemma \ref{phinapprox}, it suffices now to establish the following

\begin{proposition}
\label{vnperiodic}
Let $V$ be a finite-dimensional complex linear space,
let $S:V\to V$ be a linear operator and let $V^+\subset V$
be the subspace spanned by vectors corresponding to Jordan cells of
$S$ with eigenvalues exceeding $1$ in absolute value.
There exists a constant $C>0$ depending only on $S$ such that the following is true.
Assume that the vectors $v(n)\in V$, $n\in {\mathbb N}$, satisfy
$$
|Sv(n)-v(n+1)|<\delta
$$
for all  $n\in {\mathbb N}$ and  some constant $\delta>0$.
Then there exists a unique $v\in V^+$
such that for all $n\in {\mathbb N}$ we have
\begin{equation}
\label{slowdev}
|S^nv-v(n)|\leq {C}\cdot \delta  \cdot n^{\dim V-\dim V^+ +1}.
\end{equation}
\end{proposition}
Proof of Proposition \ref{vnperiodic}.
By definition, the subspace
$V^+$ is $S$-invariant and $S$ is invertible on $V^+$; we have furthermore that
$|Q^{-n}v|\to 0$ exponentially fast
as $n\to\infty$. Let $V^-$ be the subspace spanned by
Jordan cells corresponding to eigenvalues of absolute value at
most $1$; for $v\in V^-$, we have $|Q^nv|<Cn^{\dim V-\dim V^+}$  as $n\to\infty$.
We have the decomposition
$
V=V^+\oplus V^-.  
$
Let $$u(0)=v(0), u(n+1)=v(n+1)-Sv(n).$$
Decompose $u(n)=u^+(n)+u^-(n)$, where $u^+(n)\in V^+$,
$u^-(n)\in V^-$. Denote
$$v^+(n+1)=u^+(n+1)+Su^+(n)+\dots +S^nu^+(1);$$
$$v^-(n+1)=u^-(n+1)+Su^-(n)+\dots+S^nu^-(1);$$
$$v=u^+(0)+S^{-1}u^+(1)+\dots+S^{-n}u^+(n)+\dots.$$
By definition, $|v^-(n+1)|$ is bounded above by $C\delta n^{{\dim} V-\dim V^+ +1}$ and there exists ${\tilde C}$ 
such
that
$|S^nv-v^+(n)|<{\tilde C}\delta$
for all $n\in {\mathbb N}$, whence (\ref{slowdev}) follows. Uniqueness of $v$ follows from
the fact that for any nonzero $v^{\prime}\in V^+$ the 
sequence $|S^nv^{\prime}|$ grows exponentially as $n\to\infty$.
Proposition \ref{vnperiodic} and Lemmas \ref{thetapprox}, \ref{phinapprox} are proved completely.

Let $f\in Lip_w^+(X)$. We then have a measure $\Theta_f$  on the semi-ring ${\mathfrak C}^+_0$ given, for
$\gamma\in {\mathfrak C}^+_0$, by the formula
$$
\Theta_f(\gamma)=\int_{\gamma} fd\Phi_1^+.
$$
By (\ref{weaklip}), the measure $\Theta_f$ satisfies the assumptions of Lemma \ref{thetapprox}.
Let $\Xi^+_f\in {\cal Y}^+$ be the measure assigned to  $\Theta_f$ by Lemma \ref{thetapprox}. 

\begin{lemma}
Let  $f\in Lip_w^+(X)$,  $\Phi^-\in {\cal Y}^-_{\Gamma}$. Then
\begin{equation}
<\Xi_f^+, \Phi^->=\int_X f dm_{\Phi^-}.
\end{equation}
\end{lemma}

Proof: Choose the points  $x^{(n)}_i\in X$ satisfying (\ref{xin}).  As above,
let  ${\tilde v}\in E^-$ be the vector corresponding to $\Phi^-$ by (\ref{defPhiminusv}) and let  
${\tilde v}^{(n)}=(Q^t)^n{\tilde v}$, $n\in {\mathbb Z}$.
For any $\varepsilon>0$ and $n>0$ sufficiently large, by definition, we have
\begin{equation}
\big|m_{\Phi^-}(f)-\sum\limits_{i=1}^m \big( \int_{\gamma_n^+(x_i^{(n)})} fd\Phi_1^+  \big)\cdot
\big({\tilde v}^{(-n)}\big)_i\big|<\varepsilon.
\end{equation}

By definition of $\Xi_f^+$ and Lemma \ref{thetapprox} we have
$$
\big|\sum\limits_{i=1}^m \big( \int_{\gamma_n^+(x_i^{(n)})} fd\Phi_1^+  \big)\cdot  
\big({\tilde v}^{(-n)}\big)_i-\sum\limits_{i=1}^m \big( \Xi_f^+(\gamma_n^+(x_i^{(n)}) \big)\cdot
\big({\tilde v}^{(-n)}\big)_i|<C_{\Gamma} \cdot n^{m+1}|{\tilde v}^{(-n)}_i\big|,
$$
and, by (\ref{tildevexpfast}), the right-hand side tends to $0$ exponentially fast as $n\to\infty$.

It remains to notice that, by definition,
$$
\sum\limits_{i=1}^m \big( \Xi_f^+(\gamma_n^+(x_i^{(n)}) \big)\cdot
\big({\tilde v}^{(-n)}\big)_i=<\Xi_f^+, \Phi^->,
$$
and the Lemma is proved completely.

We have thus established that $\Xi^+_f=\Phi^+_f$, where $\Phi_f^+$ is given by (\ref{phifplus}).

\subsection{Orderings.}

Following S.~Ito \cite{ito}, A.M.~Vershik \cite{Vershik1, Vershik2}, 
assume that a partial order $\oo$ is given on ${\cal E}(\Gamma)$ in such a way that
edges starting at a given vertex are ordered linearly, while edges starting at different
vertices are not comparable. An edge will be called {\it maximal} (with respect to $\oo$)
if there does not exist a greater edge; {\it minimal}, if there does not exist a smaller edge;
and an edge $e$ will be called {\it the successor} of $e^{\prime}$
if $e>e^{\prime}$ but there does not exist $e^{\prime\prime}$ such that
$e>e^{\prime\prime}>e^{\prime}$.   

The ordering $\oo$ is extended to a partial ordering of $X_{\Gamma}$:
we write $x<x^{\prime}$ if there exists
$l\in {\mathbb Z}$ such that $x_l<x^{\prime}_l$ and $x_n=x^{\prime}_n$ for all $n>l$.
Under this ordering each leaf $\gamma^+_{\infty}$ of the foliation
${\cal F}^+$ is linearly ordered, while points lying on different leaves are not comparable.

Let $Max({\oo})$ be the set of points $x\in X$, $x=(x_n)_{n\in {\mathbb Z}}$, such that
each $x_n$ is a maximal edge.
Similarly,  $Min({\oo})$ denotes the set of points $x\in X$, $x=(x_n)_{n\in {\mathbb Z}}$, such that
each $x_n$ is a minimal edge. Since edges starting at a given vertex are ordered linearly,
the cardinalities of $Max({\oo})$ and $Min({\oo})$ do not exceed $m$.

If a leaf $\gamma^+_{\infty}$ does not intersect  $Max({\oo})$, then it does not have a maximal
element; similarly, if $\gamma^+_{\infty}$ does not intersect  $Min({\oo})$, then it does not have a
minimal element.

For $x(1), x(2)\in \gamma_{\infty}^+$, let
$$
[x(1), x(2)]= \{x^{\prime}\in \gamma_{\infty}^+: x(1)\leq x^{\prime}\leq x(2)\}.
$$

The sets $(x(1), x(2)]$, $[x(1), x(2))$, $(x(1), x(2))$ are defined similarly.

\begin{proposition}
\label{maxt}
Let $x\in X$.  If $\gamma^+_{\infty}(x)\cap Max(\oo)=\emptyset$, then for any $t\geq 0$ there exists a point $x^{\prime}\in 
\gamma^+_{\infty}(x)$ such that
\begin{equation}
\label{hplust}
\Phi_1^+([x, x^{\prime}])=t.
\end{equation}
\end{proposition}

Proof. Let $V(x)=\{t: \exists x^{\prime}\geq x: \Phi_1^+([x, x^{\prime}])=t\}$.
Since  $\gamma^+_{\infty}(x)\cap Max(\oo)=\emptyset$, for any $n$ there exists
$x^{\prime\prime}\in \gamma^+_{\infty}(x)$ such that all points in $\gamma_n^+(x^{\prime\prime})$ are greater
than $x$. Since $\Phi_1^+(\gamma_n^+(x^{\prime\prime}))$ grows exponentially, uniformly in $x^{\prime\prime}$,
as $n\to\infty$,
the set $V(x)$ is unbounded. Furthermore, since $\Phi_1^+(\gamma_n^+(x^{\prime\prime}))$ decays exponentially, uniformly in
$x^{\prime\prime}$, as $n\to-\infty$, the set $V(x)$ is dense in ${\mathbb R}_+$. Finally, by compactness of $X$, the
set $V(x)$ is closed, which concludes the proof of the Proposition.

A similar proposition, proved in the same way, holds for negative $t$. 
\begin{proposition}
\label{mint}
Let $x\in X$.  If $\gamma^+_{\infty}(x)\cap Min(\oo)=\emptyset$, 
then for any $t\geq 0$ there exists a point $x^{\prime}\in \gamma^+_{\infty}(x)$
such that
\begin{equation}
\label{hminust}
\Phi_1^+([x^{\prime}, x])=t.
\end{equation}
\end{proposition}

Define an equivalence relation $\sim$ on $X$ by writing $x\sim x^{\prime}$ if $x\in \gamma_{\infty}^+(x^{\prime})$ 
and $\Phi_1^+([x, x^{\prime}])=\Phi_1^+([x^{\prime}, x])=0$. The equivalence classes admit the following explicit description, which is 
clear from the definitions.
\begin{proposition}
\label{eqrel}
Let  $x, x^{\prime}\in X$ be such 
that $x\in \gamma_{\infty}^+(x^{\prime})$, $x<x^{\prime}$ and 
$\Phi_1^+([x, x^{\prime}])=0$.
Then there exists $n\in {\mathbb Z}$ such that 
\begin{enumerate}
\item $x^{\prime}_n$ is a successor of $x_n$; 
\item $x$ is the maximal element in $\gamma_n(x)$; 
\item  $x^{\prime}$ is the minimal element in $\gamma_n(x^{\prime})$.
\end{enumerate}
\end{proposition}
In other words, $\Phi_1^+([x, x^{\prime}])=0$ if and only if 
$(x, x^{\prime})=\emptyset$.
In particular, equivalence classes consist at most of two points and, 
$\nu$-almost surely, of only one point.

Denote $X_{\oo}=X/{\sim}$, let $\pi_{\oo}: X\to X_{\oo}$ be the projection map and set 
$\nu_{\oo}=(\pi_{\oo})_*\nu$. The probability spaces $(X_{\oo}, \nu_{\oo})$ and $(X,\nu)$ are measurably isomorphic; 
in what follows, we shall often omit the index $\oo$. The foliations ${\cal F}^+$ and ${\cal F}^-$ descend to the space $X_{\oo}$; 
we shall denote their images on $X_{\oo}$ by the same letters and, as before, denote by 
$\gamma^+_{\infty}(x)$, $\gamma_{\infty}^-(x)$ the leaves containing $x\in X_{\oo}$.  

Now let $x\in X_{\oo}$ satisfy $\gamma^+_{\infty}(x)\cap Max(\oo)=\emptyset$.
By Proposition \ref{maxt}, for any $t\geq 0$ there exists a unique $x^{\prime}$ satisfying 
(\ref{hplust}). Denote $h_t^+(x)=x^{\prime}$. 
Similarly, if $x\in X_{\oo}$ satisfy $\gamma^+_{\infty}(x)\cap Min(\oo)=\emptyset$.
By Proposition \ref{mint}, for any $t\geq 0$ there exists a unique $x^{\prime}$ satisfying
(\ref{hminust}). Denote $h_{-t}^+(x)=x^{\prime}$.

We thus obtain a flow $h_t^+$, which is well-defined on the set 
$$
X_{\oo}\setminus \Big(\bigcup\limits_{x\in Max(\oo)\cup Min(\oo)}\gamma_{\infty}^+(x)\Big),  
$$
and, in particular, $\nu$-almost surely on $X_{\oo}$. 
By (\ref{condmeas}), the flow $h_t^+$ preserves the measure $\nu$. 

More generally, it is clear from the definitions that for any $\Phi^-\in {\cal Y}^-$, 
the measure $m_{\Phi^-}$, defined by (\ref{mphi}),  
satisfies $$(h_t^+)_*m_{\Phi^-}=m_{\Phi^-},$$
similarly to G.~Forni's invariant distributions \cite{F1}, \cite{F2}.

{\bf Remark.} S.Ito in \cite{ito}
gives a construction of a flow similar to the one above. 
The flow $h_t^+$ is a continuous-time analogue of a Vershik automorphism \cite{Vershik1} 
(of which a variant also occurs in Ito's work \cite{ito}), and, in fact, is a suspension 
flow over the corresponding Vershik's automorphism, a point of view adopted in \cite{bufetov}.

\subsection{Decomposition of Arcs.}

We assume that an ordering $\oo$ is fixed on $\Gamma$.
Denote by ${\mathfrak C}(\oo)$ the semi-ring of
subsets of $X_{\Gamma}$ of the form $[x, x^{\prime})$, where $x<x^{\prime}$. 
Any measure $\Phi^+\in {\cal Y}^+$ can be extended to  ${\mathfrak C}(\oo)$ in the following way.

Let ${\mathfrak R}_n^+$ be the ring
generated
by the semi-ring ${\mathfrak C}_n^+$.
For $\gamma\in {\mathfrak C}(\oo)$, denote by $\gamma(n)$ the smallest (by inclusion) element of the ring
${\mathfrak R}_{-n}^+$ containing $\gamma$ and let ${\hat \gamma}(n)$ be the greatest (by inclusion)
element of the ring ${\mathfrak R}_{-n}^+$ contained in $\gamma$ (possibly, ${\hat \gamma}(n)=\emptyset$).
By definition,
$$
{\hat \gamma}(n)\subset {\hat \gamma}(n+1)\subset \gamma(n+1)\subset \gamma(n);
$$
\begin{equation}
\label{gammahatgamma}
\gamma(n)\setminus {\hat \gamma}(n)=\bigsqcup\limits_{i=1}^{l_n} \gamma_i^{(n)},
\end{equation}
where  $\gamma_i^{(n)}\in {\mathfrak C}_{-n}^+$, $l_n\leq ||Q||$, and 
\begin{equation}
\gamma(n)\setminus\gamma(n+1)= \bigsqcup\limits_{i=1}^{L_n} \gamma_i^{(n+1)},
\end{equation}
where $\gamma_i^{(n+1)}\in {\mathfrak C}_{-n-1}^+$, $L_n\leq 2||Q||$.

By definition, if $\Phi^+\in {\cal Y}^+$, then there are only $m$ possible values of $\Phi^+(\gamma)$ for
$\gamma\in {\mathfrak C}_{-n}^+$, and the maximum of these decays exponentially as $n\to\infty$.
We thus have

\begin{proposition}
\label{arcexpdecay}

There exists positive constants $C_{\Gamma}$, depending only on $\Gamma$, such that the following is true.
Let $v_0=0$, $v_1, \dots, v_l\in E^+$, $Qv_i=\exp(\theta)v_i+v_{i-1}$. Assume
$v\in {\mathbb C}v_1\oplus\dots\oplus {\mathbb C}v_l$ satisfies
$|v|=1$ and let $\Phi^+_v={\cal I}_{\Gamma}(v)$.
Then for any $\gamma\in {\mathfrak C}(\oo)$ we have 
$$
|\Phi^+_v(\gamma(n))-\Phi^+_v(\gamma(n+1))|\leq  C_{\Gamma}n^{l-1}\exp(-{(\Re\theta})n);
$$
$$
|\Phi^+_v({\hat \gamma}(n))-\Phi^+_v({\hat \gamma}(n+1))|\leq C_{\Gamma}n^{l-1}\exp(-({\Re\theta})n).
$$
decay exponentially as $n\to\infty$.
In particular, if $v\in E^+$, $Qv=\exp(\theta)v$, $|v|=1$, then 
$$
|\Phi^+_v(\gamma(n))-\Phi^+_v(\gamma(n+1))|\leq  C_{\Gamma}\exp(-({\Re \theta})n);
$$
$$
|\Phi^+_v({\hat \gamma}(n))-\Phi^+_v({\hat \gamma}(n+1))|\leq C_{\Gamma}\exp(-({\Re \theta})n).
$$
\end{proposition}

Consequently,  for any $\Phi^+\in {\cal Y}^+$,  $\gamma\in {\mathfrak C}(\oo)$,
the sequence $\Phi^+(\gamma(n))$ converges as $n\to\infty$, and we set
$$
\Phi^+(\gamma)=\lim\limits_{n\to\infty} \Phi^+(\gamma(n)).
$$
By (\ref{gammahatgamma}), we also have
$$
\Phi^+(\gamma)=\lim\limits_{n\to\infty} \Phi^+({\hat \gamma}(n)).
$$
\begin{proposition}
The measure $\Phi^+$ is finitely-additive on  ${\mathfrak C}(\oo)$.
\end{proposition}

Proof: Let $v\in E^+$ be such that $\Phi^+=\Phi^+_v$ and
 let $\gamma_0, \gamma_1, \dots, \gamma_k\in {\mathfrak C}(\oo)$ satisfy
$$
\gamma_0=\bigsqcup\limits_{i=1}^k \gamma_i.
$$
Consider the arcs $\gamma_0(n), \gamma_1(n), \dots, \gamma_k(n)$. We have
\begin{equation}   
\label{gammanought}
\gamma_0(n)\subset \bigcup_{i=1}^k \gamma_i(n).
\end{equation}
and decompose
$$
\gamma_i(n)=\bigsqcup \gamma_{ij}(n+1),
$$
where $\gamma_{ij}(n+1)\in  {\mathfrak C}_{-n-1}^+$.

By (\ref{gammanought}), each of the arcs $\gamma_{0j}(n+1)$ is also encountered among the arcs $\gamma_{ij}(n+1)$
(possibly, more than once, but not more than $k$ times).  Consider the collection $\gamma_{ij}(n+1)$
and cross out all the arcs  $\gamma_{0j}(n+1)$; by maximality, and since our ordering is linear on each 
leaf of the 
foliation ${\cal F}^+$, there will remain not more than  $2k||Q||$ arcs, whence we obtain
$$
\big|\sum_{i=1}^k \Phi^+(\gamma_i(n))-\Phi^+(\gamma_0(n))\big|\leq 2k||Q||\cdot |Q^{-n-1}v|,
$$
and, since the right-hand side decays exponentially as $n\to\infty$, the Proposition is proved.

\begin{lemma}
\label{mainlemma}
There exists a constant $C_{\Gamma}$ depending only on $\Gamma$ such that the following is true.
Let $f\in Lip_w^+(X_{\Gamma})$ and let $\Phi_f^+\in {\cal Y}^+$ be given by (\ref{phifplus}).
For any $\gamma\in {\mathfrak C}(\oo)$ we have
\begin{equation}
\big|\int_{\gamma} fd\Phi_1^+ - \Phi_f^+(\gamma)\big|\leq 
C_{\Gamma}||f||_{Lip_w^+}(1+\log(1+\Phi_1^+(\gamma))^{m+1}.
\end{equation}
\end{lemma}

Indeed, for $\gamma\in {\mathfrak C}^+$ this follows from Lemma \ref{thetapprox}, and for all other arcs
from Proposition \ref{arcexpdecay}.

\subsection {Ergodic averages of the flow $h_t^+$.}

Let $\Phi^+\in {\cal Y}^+$ and denote $\Phi^+[x,t]=\Phi_i^+([x, h_t^+x])$. 
The function $\Phi^+(x,t)$ is an additive cocycle over the flow $h_t^+$.
Let $f\in Lip_w^+(X_{\Gamma})$, and let $\Phi_f^+$ be defined by (\ref{phifplus}). 
By definition, $\Phi_{f\circ h_t^+}=\Phi_f^+$; recall from (\ref{invarphiplus}) that 
$\Phi_{f\circ\sigma}^+=\sigma^*\Phi_f^+$.   Lemma \ref{mainlemma} implies
\begin{theorem}
\label{multiplic}
There exists a positive constant $C_{\Gamma}$ depending only on $\Gamma$ such that
for any $f\in Lip_w^+(X_{\Gamma})$, for all $x\in X$ and all $T>0$ 
we have
$$
\big|\int_0^T f\circ h_t^+(x) dt - \Phi_f^+(x,t)\big|\leq C_{\Gamma}||f||_{Lip}(1+\log(1+T))^{m+1}.
$$
\end{theorem}

Given a bounded measurable function $f: X\to {\mathbb R}$ and $x\in X$,
introduce a continuous function ${\mathfrak S}_n[f,x]$ on the unit interval by the formula
\begin{equation}
\label{normsums}
{\mathfrak S}_n[f,x](\tau)=\int\limits_0^{\tau\exp(n\theta_1)} f\circ h^+_t(x)dt.
\end{equation}
The functions ${\mathfrak S}_n[f,x]$ are $C[0,1]$-valued random variable on the probability space $(X_{\Gamma},
\nu_{\Gamma})$.

\begin{theorem}
\label{limthm}
If $Q$ has a simple real second eigenvalue $\exp(\theta_2)$,
$\theta_2>0$, then there exists a continuous functional $\alpha: Lip_w^+(X)\to {\mathbb R}$ and a
compactly supported non-degenerate
measure $\eta$ on $C[0,1]$ such that for any $f\in Lip_{w,0}^+(X)$ satisfying $\alpha(f)\neq 0$
the sequence of random variables
$$
\frac{{\mathfrak S}_n[f,x]}{\alpha(f)\exp(n\theta_2)}
$$
converges in distribution to $\eta$ as $n\to\infty$.
\end{theorem}
  
{\bf Remark.} Compactness of the support of $\eta$ is understood in the sense of the Tchebycheff
topology on $C[0,1]$.
Nondegeneracy of the measure $\eta$  means that if $\varphi\in C[0,1]$ is distributed
according to $\eta$, then for any $t_0\in (0, 1]$ the distribution
of the real-valued random variable
$\varphi(t_0)$ is not concentrated at a single point.

The measure $\eta$ is constructed as follows: let $v_2$ be an eigenvector with eigenvalue $\exp(\theta_2)$, 
set $\Phi_2^+={\cal I}(v_2)$ (see (\ref{cali}));
then $\eta$ is the distribution of $\Phi_2^+(x, \tau)$, $0\leq \tau\leq 1$, considered as a $C[0,1]$-valued 
random variable on the space $X_{\Gamma}, \nu_{\Gamma})$.
The functional $\alpha(f)$ is constructed as follows: under the assumptions of Theorem \ref{limthm}, 
the matrix $Q^t$ also has the simple real second eigenvalue $\exp(\theta_2)$; let ${\tilde v}_2$ be the  
eigenvector with eigenvalue $\exp(\theta_2)$, normalized in such a way that $\sum_{i=1}^m (v_2)_i ({\tilde v}_2)_i=1$; 
set $\Phi_2^-={\tilde {\cal I}}({\tilde v}_2)$ (see (\ref{tildecali})),and let $m_{\Phi_2^-}$ be given by (\ref{mphi}); then  
$$
\alpha(f)=\int fdm_{\Phi_2^-}.
$$

\subsection{The diagonalizable case.}

As an illustration, consider the case when $Q|_{E^+}$ 
is diagonalizable with eigenvalues $\exp({\theta_i})$, $i=1, \dots, r$,
$\Re(\theta_i)>0$. The Perron-Frobenius vector $h$ corresponds to $\exp(\theta_1)$; 
let $v_2, \dots, v_r$ be eigenvectors corresponding to $\exp(\theta_i)$: thus  $Qv_i=\exp(\theta_i)v_i$, $i=2, \dots, r$ and
$$
E^+={\mathbb C}h\oplus {\mathbb C}v_2\oplus \dots\oplus {\mathbb C}v_r
$$
We have a similar direct-sum representation  for $Q^t$:
$$
{\tilde E}^+={\mathbb C}\la\oplus {\mathbb C}{\tilde v}_2\oplus \dots\oplus 
{\mathbb C}{\tilde v}_r,
$$
where  $Q^t{\tilde v}_i=\exp(\theta_i){\tilde v}_i$,  $i=2, \dots, r$.
For $i\neq j$ we have 
\begin{equation}
\sum_{l=1}^m ({v}_i)_l({\tilde v}_j)_l=0,
\end{equation}
and, for normalization, let us assume that for all $i=1, \dots, r$ we have  
\begin{equation}
\sum_{l=1}^m ({v}_i)_l({\tilde v}_i)_l=1.       
\end{equation}

Let $\Phi^+_i={\cal I}(v_i)$, $\Phi^-_i={\tilde {\cal I}}({\tilde v}_i)$, $i=2, \dots, r$.
Since $\Phi_1^+={\cal I}(h)$, the measures $\Phi^+_i$,  $i=1, \dots, r$, form a basis in 
${\cal Y}^+$, for which the measures $\Phi_1^-={\tilde {\cal I}}(\la)$, $\Phi^-_2, \dots, \Phi^-_r$ 
form a dual basis in ${\cal Y}^-$. 

For $i=1, \dots, r$, from (\ref{mphi}) we have the measures 
$m_{\Phi^-_i}=\Phi_1^+\times \Phi_i^-$. For instance, $m_{\Phi_1^-}=\nu$. 
Theorem \ref{multiplic} now implies  
\begin{corollary}
\label{diagmultiplic}
For any $f\in Lip_w^+(X_{\Gamma})$ we have 
$$
\big|\int_0^T f\circ h_t^+(x) dt - T\int_X fd\nu - \sum_{i=2}^r 
\Phi_i^+(x,T)\big(m_{\Phi_i^-}(f)\big)\big|\leq 
C_{\Gamma}||f||_{Lip}(1+\log(1+T))^{m+1},   
$$
where $C_{\Gamma}$ is a constant depending only on $\Gamma$.
\end{corollary}

For the action of the shift we have: 
\begin{equation}
\label{sigmaplus}
(\sigma)_*\Phi^+_i=\exp(-\theta_i)\Phi^+_i, \ i=1, \dots, r; 
\end{equation}
\begin{equation}
\label{sigmaminus}
(\sigma)_*\Phi^-_i=\exp(\theta_i)\Phi^-_i, \ i=1, \dots, r. 
\end{equation}
Corollary \ref{diagmultiplic} now yields 
\begin{equation}
\label{sigman}
\int\limits_0^{\tau\exp(\theta_1n)} f\circ h_t^+(x) dt=\sum_{i=1}^r \exp(n\theta_i) 
m_{\Phi_i^-}(f)\Phi_i^+(\sigma^nx, \tau)+O(n^{m+1}). 
\end{equation}

\subsection{The H{\"o}lder property.}

As above, we write $\Phi^+(x,t)=\Phi^+([x, h_t^+x])$. Our next aim is to show that $\Phi^+(x,t)$ is H{\"o}lder 
in $t$ for any $x\in X_{\oo}$.
\begin{proposition}
\label{hoeldereigen}
There exist positive constants $C_{\Gamma}$ and $t_0$, depending only on $\Gamma$ such that the following is true.
Let $v\in E^+$, $Qv=\exp(\theta)v$, $|v|=1$. Then for all $x\in X$ and positive $t<t_0$ we have 
$$
|\Phi_v^+(x,t)|\leq C_{\Gamma}t^{\Re\theta/\theta_1}.
$$ 
\end{proposition}
\begin{proposition}
\label{hoelderjordan}
There exist positive constants $C_{\Gamma}$ and $t_0$, depending only on $\Gamma$ such that the following is true.
Let $v_0=0$, $v_1, \dots, v_l\in E^+$, $Qv_i=\exp(\theta)v_i+v_{i-1}$. Assume 
$v\in {\mathbb C}v_1\oplus\dots\oplus {\mathbb C}v_l$ satisfies
$|v|=1$. Then for all $x\in X$ and positive $t<t_0$ we have 
$$
|\Phi_v^+(x,t)|\leq C_{\Gamma}|\log t|^{l-1}t^{\Re\theta/\theta_1}.
$$
\end{proposition}

Proof of Propositions \ref{hoeldereigen}, \ref{hoelderjordan}. Denote $\gamma=[x, h^+_tx]$. If $t$ is small enough, then 
${\hat \gamma}(0)=\emptyset$.
Let $n_0$ be the smallest positive integer such that ${\hat \gamma}(n_0)\neq\emptyset$. There exist  
positive constants $C_1, C_2$, depending only on $\Gamma$, such that
$$
C_1 t\leq \exp(-\theta_1n_0)\leq C_2 t,
$$
and Propositions  \ref{hoeldereigen}, \ref{hoelderjordan} follow now from Proposition \ref{arcexpdecay}.

\begin{corollary}
\label{hoelderuniform}
There exist positive constants $\theta>0$ and $t_0>0$ depending only on $Q$ such that 
for all $v\in E^+$, $|v|=1$, all $x\in X$ and all positive $t<t_0$ we have 
$$
|\Phi_v^+(x,t)|\leq t^{\theta/\theta_1}.
$$ 
\end{corollary}

For $v\in E^+$, $|v|=1$ denote
$$
\theta_v=\lim_{n\to\infty}\frac{\log|Q^nv|}{n}.
$$
\begin{corollary}
\label{upperphi}
For any $\varepsilon>0$ there exists a constant $T_{\varepsilon}$ depending only on $\varepsilon$ and $\Gamma$ 
such that for any $v\in E^+$, $|v|=1$, any $x\in X$ and any $T>T_{\varepsilon}$, we have 
$$
|\Phi^+_v(x,T)|\leq T^{\theta_v/\theta_1+\varepsilon}.
$$
\end{corollary}

Proof: Indeed, let $t_0$ be the constant given by Proposition\ref{hoelderjordan}. Let $n_0=n_0(T)$ 
be the smallest such integer that  $T=\tau\exp(n(T)\theta_1)$, where $\tau<t_0$. Since 
$\Phi^+_v(x,T)=\Phi^+_{Q^nv}(\sigma^nx, \tau)$ for all $n$, 
it follows from Proposition \ref{hoelderjordan} that 
$$
|\Phi^+(x,T)|\leq C_{\Gamma} n_0^{m+1}\exp(n_0\Re(\theta_v))\leq C_{\Gamma}  T^{\theta_v/\theta_1+\varepsilon}
$$
if $T$ is sufficiently large (depending only on $\varepsilon$).

\begin{corollary}
\label{logphi}
For any $v\in E^+$ we have
\begin{equation}
\label{limsupphiv}
\limsup_{T\to\infty} \frac{\log |\Phi^+_v(x,T)|}{\log T}=\frac{\theta_v}{\theta_1}.
\end{equation}
\end{corollary}
Indeed, the upper bound for the limit superior follows from Corollary \ref{upperphi}, and the lower bound
is immediate from the relation $\Phi^+_v(\gamma_n(x))=(Q^nv)_{F(x_{n+1})}$.

\begin{corollary}
\label{nonconst}
For any $\tau\in {\mathbb R}$ and any $v\in E^+$ satisfying $v\neq 0$, $\sum\limits_{i=1}^m v_i\la_i=0$, the function 
$\Phi_v^+(x,\tau)$ 
is not a constant in $x$.
\end{corollary}
Proof: Indeed, assume $\Phi_v^+(x,\tau)=c$ identically. Then $\Phi^+(x, k\tau)=kc$, which contradicts  
(\ref{limsupphiv}): is $c=0$, then the limit superior is $0$; if $c\neq 0$, then the 
limit superior is $1$.

\subsection{Tightness.}

In this subsection, we assume that $Q$ has a simple real second eigenvalue $\exp(\theta_2)$,
$\theta_2>0$. Let $v_2$ be the corresponding eigenvector and let $\Phi_2^+={\cal I}(v_2)$.
Take $x\in X$ and consider $\Phi^+(x,\tau)$ as a continuous function of $\tau$ on the unit interval.
Let $\eta$ be the distribution of $\Phi_2^+(x,\tau)$ in $C[0,1]$.
Note that by Corollary \ref{nonconst}, for any $\tau_0$ the value of $\Phi_2^+(x,\tau)$ is not constant on $X$, 
so the measure $\eta$ is nondegenerate.

Let ${\mathfrak S}_n[f, x]$ be defined by the equation (\ref{normsums}).
Introduce a sequence of measures $\mu_n$ on $C[0,1]$ by the formula
$\mu_n={\mathfrak S}[n,f]_*\nu_{\Gamma}$.

By Theorem 8.1 in Billingsley \cite{billingsley}, p.54, to prove Theorem \ref{limthm} it 
suffices to establish the following two Lemmas. 
\begin{lemma}
\label{findim}
Finite-dimensional distributions of the measures $\mu_n$ weakly converge to those of $\eta$.
\end{lemma}
 \begin{lemma}
\label{tight}
The family $\mu_n$ is tight in $C[0,1]$.
\end{lemma}
Proof of Lemma \ref{findim}. By Theorem \ref{multiplic}
$$
\int_0^T f\circ h_t^+(x)dt=\Phi_f^+(x,T)+O((\log T)^{m+1}).
$$
Let $v_2$ be the eigenvector corresponding to the eigenvalue $\exp(\theta_2)$, $|v|=1$,
and let $\Phi_2^+\in {\cal Y}^+$ be the corresponding measure.
We have
$$
E^+={\mathbb C}v_2\oplus E_3,
$$
where $E_3$ is spanned by Jordan cells corresponding to eigenvalues with absolute value
less than $\exp(\theta_2)$. Let $\zeta$ be a number smaller than $\theta_2$ but greater than
the spectral radius of $Q|_{E_3}$.
Write
\begin{equation}
\label{alphabeta}
\Phi_f^+=\alpha(f)\Phi_2^+ + \beta(f) \Phi^+_{v_3},
\end{equation}
where $v_3\in E^+$, $|v_3|=1$, and $\alpha(f), \beta(f)$ are continuous functionals on  $Lip_w^+(X)$, so, in
particular, we have
$$
|\alpha(f)|<C_{01}||f||_{Lip_w^+}; \ |\beta(f)|<C_{02}||f||_{Lip_w^+},
$$
where the constants $C_{01}$, $C_{02}$ only depend on $\Gamma$.

By Corollary \ref{hoelderuniform}, there exists $t_0$ depending only on $\Gamma$ such that for
any positive $t$ such that $t<t_0$, any $x\in X$ and any $v\in E^+$ satisfying $|v|=1$ we have
\begin{equation}
\label{leqone}
|\Phi_v^+(x,t)|\leq 1.
\end{equation}
  
Write $T=t\exp(n\theta_1)$, where $t<t_0$.
Since $\Phi_{v_3}^+(x,T)=\Phi^+_{Q^nv_3}(\sigma^nx, t)$,
for all sufficiently large $n$, we have $|Q^nv_3|<\exp(\zeta n)$ and therefore
\begin{equation}
\label{phi3small}
|\Phi_{v_3}^+(x,\tau \exp(n\theta_1))|< \exp(n\zeta)
\end{equation}
for all $x\in X$.    
By Theorem \ref{multiplic} we have 
\begin{equation}
\big|\int\limits_0^{\tau \exp(n\theta_1)} f\circ h_t^+(x)dt - 
\Phi_f^+(x,\tau\exp(\theta_1n))\big|=O(n^{m+1}). 
\end{equation}
Since 
$$
\Phi_f^+(x, \tau\exp(n\theta_1))=\alpha(f)\Phi_2^+((x, \tau\exp(n\theta_1)) + 
\beta(f) \Phi^+_{v_3}(x, \tau\exp(n\theta_1))
$$
combining the equality 
$$
\Phi_2^+(x, \tau\exp(n\theta_1))=\exp(n\theta_2)\Phi_2^+(\sigma^nx,\tau) 
$$
with the bound (\ref{phi3small}), we obtain, for all large $n$ and all $x\in X$, uniformly in $\tau\in [0,1]$,  
the estimate
$$
|{\mathfrak S}_n[f,x](\tau)- \alpha(f) \Phi_2^+(\sigma^nx,\tau)|\leq C_{\Gamma} 
||f||_{Lip_w^+}\exp((\zeta-\theta_2)n).
$$
Since $\sigma$ preserves the measure $\nu$, it follows that the $k$-dimensional distributions of 
$\big({\mathfrak S}_n[f,x](\tau_1), {\mathfrak S}_n[f,x](\tau_2), \dots, {\mathfrak S}_n[f,x](\tau_k)\big)$
converge to the $k$-dimensional distribution of 
$\big(\Phi_2^+(x,\tau_1), \Phi_2^+(x,\tau_2), \dots, \Phi_2^+(x,\tau_k)\big)$, and 
Lemma \ref{findim} is proved.

The argument above yields also 
\begin{proposition}
\label{unifest}
There exist positive constants $C_0=C_0(\Gamma)$ and $T_0=T_0(\Gamma)$
such that for any $x\in X$, any $f\in Lip_{w,0}^+(X)$  and any $T>T_0$
we have
$$
|\int_0^T f\circ h_t^+(x)dt|\leq C_0\cdot ||f||_{Lip_w^+}\cdot T^{\theta_2/\theta_1}.
$$
\end{proposition}

Indeed,
for sufficiently large $T$,  $T=t\exp(n\theta_1)$, where $t<t_0$, 
from (\ref{alphabeta}) we have
$$
\Phi_f^+(x,T)=\alpha(f)\exp(n\theta_2)\Phi_2^+(\sigma^nx, t)+ O(\exp(n{\zeta})).
$$
Since, by (\ref{leqone}), we have $|\Phi_2^+(\sigma^nx, t)|\leq 1$, Proposition \ref{unifest} is established.

We proceed to the proof of Lemma \ref{tight}.
\begin{proposition}
\label{veryuniformhoelder} 
There exists a constant $C_{\Gamma}$ depending only on $\Gamma$ 
such that for any $f\in Lip_{w,0}^+(X)$,  
any $n>0$, any $x\in X$ and any $\tau_1, \tau_2\in [0,1]$, we have
$$
|{\mathfrak S}_n[x,f](\tau_2)-{\mathfrak S}_n[x,f](\tau_1)|\leq 
C_{\Gamma}||f||_{Lip_w^+}|\tau_2-\tau_1|^{\theta_2/\theta_1}.
$$ 
\end{proposition}

Lemma \ref{tight} follows from Proposition 
\ref{veryuniformhoelder} by the Arzel{\`a}-Ascoli Theorem.

Proof of Proposition \ref{veryuniformhoelder}:
Let $\tau_1, \tau_2\in [0,1]$, $\tau_1<\tau_2$. 
For brevity, write ${\mathfrak S}_n={\mathfrak S}_n[f,x].$
We have then
$$
{\mathfrak S}_n(\tau_2)-{\mathfrak S}_n(\tau_1)=\frac 1{\exp(n\theta_2)}\int\limits_{\tau_1\exp(n\theta_1)}^{\tau_2\exp(n\theta_1)} 
f\circ h_t^+(x)dt.
$$
Let $T_0$ be the constant given by Proposition 
\ref{unifest} and assume first that 
$$
(\tau_2-\tau_1)\cdot \exp(n\theta_1)\geq T_0.
$$ 
By Proposition \ref{unifest} we have
$$
\int\limits_{\tau_1\exp(n\theta_1)}^{\tau_2\exp(n\theta_1)} f\circ h_t(x)dt\leq
C||f||_{Lip_w^+}\cdot (\tau_2-\tau_1)^{\theta_2/\theta_1}\exp(n\theta_2),
$$
and, consequently,
$$
|{\mathfrak S}_n(\tau_2)-{\mathfrak S}_n(\tau_1)|\leq C_{33}(\tau_2-\tau_1)^{\theta_2/\theta_1},
$$
where the constant $C_{33}$ only depends on $\Gamma$.  

Now let $\tau_2-\tau_1=\tau_0\exp(-n\theta_1)$, $\tau_0<T_0$.
Since 
$$
\exp(-n\theta_2)=((\tau_2-\tau_1)/\tau_0)^{\theta_2/\theta_1},
$$ 
using boundedness of $f$, write
$$
\frac 1{\exp(n\theta_2)}\int\limits_{\tau_1\exp(n\theta_1)}^{\tau_2\exp(n\theta_1)} f\circ h_t^+(x)dt
\leq \exp(-n\theta_2)\cdot ||f||_{\infty}\cdot \tau_0\leq 
$$
$$
\leq
\tau_0^{1-\theta_2/\theta_1}||f||_{\infty}(\tau_2-\tau_1)^{\theta_2/\theta_1}\leq
T_0^{1-\theta_2/\theta_1}||f||_{\infty}(\tau_2-\tau_1)^{\theta_2/\theta_1},
$$
and the Proposition is proved.
Theorem \ref{limthm} is proved completely.

\subsection{A symbolic coding for translation flows on surfaces.}

To derive Theorems \ref{multiplicpsan}, \ref{limthmpsan} from Theorems \ref{multiplic}, \ref{limthm}, 
it remains to observe that the vertical flow on the stable foliation of a pseudo-Anosov diffeomorphism is 
isomorphic to a symbolic flow on the asymptotic foliation of a Markov compactum obtained from  the 
decomposition of the underlying surface into Veech's zippered rectangles, see \cite{bufetov}, Sec. 4.
The identification of $E^+$ (and, consequently, of ${\cal Y}^+$)
with the corresponding subspace in cohomology is given by Proposition 4.16 in Veech\cite{V3}.   
The fact that the pairing between cocycles corresponds to the cup-product is immediate 
from Proposition 4.19 in \cite{V3}. 

\section{Spaces of  Markov Compacta.}

Let ${\mathfrak G}$ be the set of all oriented graphs on $m$ vertices such that there is an edge starting at 
every vertex and an edge ending at every vertex. As before, for a graph $\Gamma\in {\mathfrak G}$, we denote 
by ${\cal E}(\Gamma)$ the set of its edges and by
$A(\Gamma)$ its incidence matrix: $A_{ij}(\Gamma)=\# \{e\in {\cal E}(\Gamma): I(e)=i, F(e)=j\}.$
Denote $\Omega={\mathfrak G}^{\mathbb Z}$:
$$
\Omega=\{\omega=\dots \omega_{-n}\dots\omega_n\dots, \omega_i\in {\mathfrak G}, i\in {\mathbb Z} \},
$$

For $\omega\in\Omega$, denote by $X(\omega)$ the corresponding Markov compactum:
$$
X(\omega)=\{x=\dots x_{-n}\dots x_n \dots, x_n\in {\cal E}(\omega_n), F(x_{n+1})=I(x_n)\}.
$$

For $x\in X$, $n\in {\mathbb Z}$, introduce the sets
$$
\gamma^+_n(x)=\{x^{\prime}\in X(\omega): x^{\prime}_t=x_t, t\geq n\};\
\gamma^-_n(x)=\{x^{\prime}\in X(\omega): x^{\prime}_t=x_t, t\leq n\};
$$
$$
\gamma^+_{\infty}(x)=\bigcup_{n\in {\mathbb Z}} \gamma_n^+(x);\
\gamma^-_{\infty}(x)=\bigcup_{n\in {\mathbb Z}} \gamma_n^-(x).
$$
The sets $\gamma^+_{\infty}(x)$ are leaves of the asymptotic foliation  ${\cal F}^+_{\omega}$
on $X(\omega)$; the sets $\gamma^-_{\infty}(x)$ 
are leaves of the asymptotic foliation  ${\cal F}^-_{\omega}$ on $X(\omega)$.

For $n\in {\mathbb Z}$ let ${\mathfrak C}^+_{n, \omega}$ be the collection of all subsets of $X(\omega)$ 
of the 
form
$\gamma^+_n(x)$, $n\in {\mathbb Z}$, $x\in X$; similarly,
${\mathfrak C}^-_{n, \omega}$ is the collection of all subsets of the form  $\gamma^-_n(x)$.
Set
\begin{equation}
\label{cplusom}
{\mathfrak C}^+_{\omega}=\bigcup\limits_{n\in {\mathbb Z}} {\mathfrak C}^+_{n, \omega}; 
{\mathfrak C}^-_{\omega}=\bigcup\limits_{n\in {\mathbb Z}} {\mathfrak C}^-_{n, \omega}.
\end{equation}

Just as in the periodic case, the collections ${\mathfrak C}^+_{n, \omega}$, 
${\mathfrak C}^-_{n, \omega}$, ${\mathfrak C}^+_{\omega}$, ${\mathfrak C}^-_{\omega}$  are  semi-rings.   

{\bf Remark.} To make notation lighter, we shall often omit the subscript $\omega$ 
and only include it when dependence on $\omega$ is underlined.

\subsection{Measures and Cocycles.}

Let $\sigma$ be the shift  on $\Omega$ given by the formula $(\sigma\omega)_n=\omega_{n+1}$. 
Let $\Prob$ be an ergodic $\sigma$-invariant probability measure on $\Omega$.
We then have a natural cocycle ${\mathbb A}$ on the system
$(\Omega, \sigma, \Prob)$
defined, for $n>0$, by the formula
$$
{\mathbb A}(n,\omega)=A(\omega_{n})\dots A(\omega_1).
$$
The cocycle ${\mathbb A}$ will be called the {\it renormalization cocycle}.

We need the following assumptions on the measure $\Prob$ and on the cocyle ${\mathbb A}$.

\begin{assumption}
\label{astwo}
The matrices $A(\omega_n)$ are almost surely invertible with respect to $\Prob$.
There exists $\Gamma \in {\mathfrak G}$ such that $\Prob(\Gamma)>0$.
\end{assumption}
\begin{assumption}
\label{asos}
The logarithm of the renormalization cocycle
(and of its inverse) is integrable.
\end{assumption}

For $n<0$ set
$$
{\mathbb A}(n,\omega)=A^{-1}(\omega_{-n})\dots A^{-1}(\omega_0).
$$
and set ${\mathbb A}(0,\omega)$ to be the identity matrix.

The {\it transpose} cocycle  ${\mathbb A}^t$ over the dynamical 
system $(\Omega, \sigma^{-1}, \Prob)$
defined, for $n>0$, by the formula
$$
{\mathbb A}^t(n,\omega)=A^t(\omega_{1-n})\dots A^t(\omega_0).
$$
Similarly, for $n<0$ write
$$
{\mathbb A}^t(n,\omega)=(A^t)^{-1}(\omega_{-n})\dots (A^t)^{-1}(\omega_1).
$$
and set ${\mathbb A}^t(0,\omega)$ to be the identity matrix.

By Assumptions \ref{astwo}, \ref{asos}, for $\Prob$-almost any $\omega\in\Omega$ we have the 
decompositions
$$
{\mathbb R}^m=E^+_{\omega}\oplus E^-_{\omega};  \ {\mathbb R}^m={\tilde E}^+_{\omega}\oplus 
{\tilde E}^-_{\omega},
$$
where $E^+$ is the Lyapunov subspace corresponding to positive Lyapunov exponents of ${\mathbb A}$; 
${\tilde E}^+$ is the Lyapunov subspace corresponding to positive Lyapunov exponents of ${\mathbb A^t}$;  
$E^-$ is the Lyapunov subspace corresponding to zero and negative Lyapunov exponents of ${\mathbb A}$;
$E^-$ is the Lyapunov subspace corresponding to zero and negative Lyapunov exponents of ${\mathbb A}^t$.
The standard inner product on ${\mathbb R}^m$ yields a nondegenerate 
pairing between the spaces $E^+_{\omega}$ and ${\tilde E}^+_{\omega}$.

In particular, by Assumption 1, the spaces $E^+_{\omega}$ and ${\tilde E}^+_{\omega}$ each contain a unique vector all 
whose coordinates are positive; we denote these vectors by $h^{(\omega)}$ and $\la^{(\omega)}$, 
respectively, and assume that they are normalized by (\ref{lahnorm}).

Let $v\in E^+_{\omega}$ and for all $n\in {\mathbb Z}$ set $v^{(n)}={\mathbb A}(n, \omega)v$.
Introduce a
finitely-additive complex-valued measure $\Phi^+_v$ on the semi-ring 
${\mathfrak C}^+_{\omega}$ (defined in (\ref{cplusom})) by the formula
\begin{equation}
\label{defPhiplusvom}
\Phi^{+}_v(\gamma_{n+1}^+(x))=(v^{(n)})_{F(x_{n+1})}.
\end{equation}

As before, the measure $\Phi^+_v$ is invariant under 
holonomy along ${\cal F}^-$: by definition, we have 
the following
\begin{proposition}
\label{holoplusom}
If $F(x_n)=F(x^{\prime}_n)$, then $\Phi^{+}_v(\gamma_n^+(x))=\Phi^{+}_v(\gamma_n^+(x^{\prime}))$.
\end{proposition}

The measures $\Phi^+_v$ span a complex linear space, which is denoted ${\cal Y}^+_{\omega}$.
The map ${\cal I}_{\omega}: v \to \Phi^+_v$ is an isomorphism between $E^+_{\omega}$ and  ${\cal 
Y}^+_{\omega}$. Set $\Phi_{1, \omega}^+={\cal I}_{\omega}(h^{(\omega)})$.

Now for ${\tilde v}\in {\tilde E}^+$ and for all $n\in {\mathbb Z}$ 
set ${\tilde v}^{(n)}={\mathbb A}^t(n, \omega){\tilde v}$   
and  introduce a
finitely-additive complex-valued measure $\Phi^-_{{\tilde v}}$ on the semi-ring ${\mathfrak 
C}^-_{\omega}$
(defined in (\ref{cplusom})) by the formula
\begin{equation}
\label{defPhiminusvom}
\Phi^{-}_{{\tilde v}}(\gamma_n^-(x))=({\tilde v}^{(-n)})_{I(x_n)}.
\end{equation}

By definition, the measure  $\Phi^-_{{\tilde v}}$ is invariant under holonomy along ${\cal F}^+$:
more precisely, we have the following
\begin{proposition}
\label{holominusom}
If $I(x_n)=I(x^{\prime}_n)$, then 
$\Phi^{-}_{{\tilde v}}(\gamma_n^-(x))=\Phi^{-}_{{\tilde v}}(\gamma_n^-(x^{\prime}))$.
\end{proposition}

Let  ${\cal Y}^-_{\omega}$ be the space spanned by the measures
$\Phi^-_{{\tilde v}}$, ${\tilde v}\in {\tilde E}^+$.
The map ${\tilde {\cal I}}_{\omega}: {\tilde v} \to \Phi^-_{{\tilde v}}$ is an isomorphism
between ${\tilde E}^+_{\omega}$ and ${\cal Y}^-_{\omega}$.  Set $\Phi_{1, \omega}^-={\tilde {\cal I}}_{\omega}(\la^{(\omega)})$.

Define a map $t_{\sigma}: X_{\omega}\to X_{\sigma\omega}$ by $(t_{\sigma} x)_i=x_{i+1}$.
The map $t_{\sigma}$ induces a map $t_{\sigma}^*:{\cal Y}^+_{\sigma\omega}\to{\cal Y}^+_{\omega}$ 
given, for $\Phi^+_{\sigma\omega}\in {\cal Y}^+_{\sigma\omega}$ and 
$\gamma\in {\mathfrak  C}^+_{\omega}$, 
by the formula 
$$
t_{\sigma}^*\Phi^+(\gamma)=\Phi^+_{\sigma\omega}(t_{\sigma}\gamma).
$$ 
We have the following commutative diagrams:
$$
\begin{CD}
E^+_{\omega}@ >{\cal I}_{\omega}>> {\cal Y}^+_{\omega} \\
@VV{\mathbb A}(1,\omega)V           @AAt_{\sigma}^*A   \\
E^+_{\sigma\omega}@ >{\cal I}_{\sigma\omega}>> {\cal Y}^+_{\sigma\omega} \\
\end{CD}
$$
$$
\begin{CD}
{\tilde E}^+_{\omega}@ >{\tilde {\cal I}}_{\omega}>> {\cal Y}^-_{\omega} \\
@AA{\mathbb A}^t(1,\sigma\omega)A           @AAt_{\sigma}^*A   \\
{\tilde E}^+_{\sigma\omega}@ >{\tilde {\cal I}}_{\sigma\omega}>> {\cal Y}^-_{\sigma\omega} \\
\end{CD}
$$
\subsection{Pairings and weakly Lipschitz functions.}
Given $\Phi^+\in {\cal Y}^+_{\omega}$, $\Phi^-\in {\cal Y}^-_{\omega}$, introduce
a finitely additive measure $\Phi^+\times \Phi^-$ on the semi-ring ${\mathfrak C}$
of cylinders in $X(\omega)$: for any $C\in {\mathfrak C}$ and $x\in  C$, set
\begin{equation}
\label{prodmeasgen}
\Phi^+\times \Phi^-(C)=\Phi^+(\gamma^+_{\infty}(x)\cap C)\cdot 
\Phi^-(\gamma^-_{\infty}(x)\cap C).
\end{equation}
Note that by Propositions \ref{holoplusom}, \ref{holominusom}, 
the right-hand side in (\ref{prodmeasgen}) does not depend on $x\in {C}$.

As above, for $\Phi^-\in {\cal Y}^-_{\omega}$, denote                  
\begin{equation}
\label{mphigen}
m_{\Phi^-}=\Phi_1^+\times \Phi^-.
\end{equation}
In particular, we have a positive countably additive measure 
$$
\nu_{\omega}=\Phi^+_{h^{(\omega)}}\times \Phi^-_{\la^{(\omega)}}.
$$

There is a natural ${\mathbb C}$-linear pairing $<,>$ between the spaces 
${\cal Y}^+_{\omega}$ and 
${\cal Y}^-_{\omega}$: for $\Phi^+\in {\cal Y}^+_{\omega}$, 
$\Phi^-\in {\cal Y}^-_{\omega}$, set
\begin{equation}
<\Phi^+, \Phi^->=\Phi^+\times \Phi^-(X(\omega)).
\end{equation}
As in Sec. \ref{secpair}, we have 
\begin{proposition}
Let  $v\in E^+_{\omega}$, ${\tilde v}\in {\tilde E}^+_{\omega}$, 
$\Phi^+_v={\cal I}_{\omega}(v)$,
$\Phi^-_{\tilde v}={\tilde {\cal I}}_{\omega}({\tilde v})$. Then
\begin{equation}
<\Phi^+_v, \Phi^-_{{\tilde v}}>=\sum_{i=1}^m v_i{\tilde v}_i.
\end{equation}
The pairing $<,>$ is non-degenerate and $t_{\sigma}^*$-invariant.
\end{proposition}
  
The function space  $Lip_w^+(X(\omega))$ is introduced in the same way as before:
a bounded Borel-measurable function $f:X(\omega) \to {\mathbb C}$ belongs to the space $Lip_w^+(X)$ 
if there exists a 
constant $C>0$ such that for all $n\geq 0$ and any $x, x^{\prime} \in X$ satisfying
$F(x_{n+1})=F(x^{\prime}_{n+1})$, we have
\begin{equation}
\label{weaklipgen} 
|\int_{\gamma_n^+(x)}fd\Phi_1^+  - \int_{\gamma_n^+(x^{\prime})} fd\Phi_1^+|\leq C,
\end{equation}
and, if $C_f$ is the infimum of all $C$ satisfying (\ref{weaklipgen}), then we norm $Lip_w^+(X)$ by 
setting
$$
||f||_{Lip_w^+}=\sup_X f+C_f.
$$
As before, we denote by $Lip_{w,0}^+(X(\omega))$ the subspace of functions of 
$\nu_{\omega}$-integral zero.

Take $\Phi^-\in {\cal Y}^-$. Any function $f\in Lip_w^+(X)$ is integrable with respect
to the measure $m_{\Phi^-}$ in the same sense as in Sec. \ref{secweaklip}, 
and a measure $\Phi_f^+\in {\cal Y}^+$ is defined by the requirement 
that for any $\Phi^-\in {\cal Y}^-$ we have
\begin{equation}
\label{phifplusgen}
<\Phi_f^+, \Phi^->=\int_{X(\omega)} f dm_{\Phi^-}.
\end{equation}

Note that the mapping ${\Xi}^+_{\omega}:  Lip_w^+(X(\omega))\to {\cal Y}^+_{\omega}$ 
given by
${\Xi}^+_{\omega}(f)=\Phi^+_f$ is continuous by definition  and satisfies
\begin{equation}
\label{invarphiplusgen} 
{\Xi}^+_{\sigma\omega}(f\circ t_{\sigma})=(t_{\sigma})^*{\Xi}^+_{\omega}(f).
\end{equation}  

From the definitions we also have
\begin{proposition}
Let $\Phi^+(1), \dots, \Phi^+(r)$ be a basis in ${\cal Y}^+_{\omega}$
and let $\Phi^-(1), \dots, \Phi^-(r)$ be the dual basis in ${\cal Y}^-_{\omega}$ with respect to 
the pairing $<,>$.
Then  for any $f\in Lip_w^+(X(\omega))$ we have
$$
\Phi_f^+=\sum_{i=1}^r \big(m_{\Phi^-(i)}(f)\big)\Phi^+(i).
$$
\end{proposition}

\subsection{Orderings and flows.}

Assume that for $\Prob$-almost every $\omega$ a partial  ordering $\oo(\omega)$ is given on 
${\cal E}(\omega_n)$ for all $n\in {\mathbb Z}$ in such a way that edges starting at a given vertex 
are ordered linearly, while edges starting at different
vertices are incomparable. Assume, moreover, that the orders $\oo(\omega)$ are $\sigma$-invariant, 
in the sense that the ordering  $\oo(\omega)$ on ${\cal E}(\omega_n)$ is the same as 
the ordering $\oo(\sigma\omega)$ on  ${\cal E}((\sigma\omega)_{n-1})$.

Similarly to the above, construct spaces $X_{\oo}(\omega)$ and introduce a flow $h_t^{(+, \omega)}$ 
on each   $X_{\oo}(\omega)$. The shift $\sigma$ renormalizes the flows $h_t^{(+, \omega)}$:
if we set
\begin{equation}
\label{hone}
H^{(1)}(n, \omega)=||{\mathbb A}(n, \omega)||,
\end{equation}
then for any $t\in {\mathbb R}$ we have a commutative diagram
$$
\begin{CD}
X(\omega)@ >h_t^{(+, \omega)}>> X(\omega) \\
@VV t_{\sigma}V           @VV t_{\sigma}V   \\
X(\sigma\omega)@ >h_{t/H^{(1)}(1,\omega)}^{(+, \sigma\omega)}>> X(\sigma\omega) \\
\end{CD}
$$

As before, each measure $\Phi^+\in {\cal Y}^+_{\omega}$ yields a H{\"o}lder cocycle over the flow 
$h_t^{(+, \omega)}$; we shall denote the cocycle by the same letter as the measure.

Note that for any $\Phi^-\in {\cal Y}^-_{\omega}$ the measure $m_{\Phi^-}$ defined by
(\ref{mphigen}) satisfies
$$(h_t^{(+,\omega)})_*m_{\Phi^-}=m_{\Phi^-},$$
similarly to G. Forni's invariant distributions \cite{F1}, \cite{F2}.

Note that the mapping ${\Xi}^+_{\omega}:  Lip_w^+(X(\omega))\to {\cal Y}^+_{\omega}$
given by ${\Xi}^+_{\omega}(f)=\Phi^+_f$ by definition satisfies
\begin{equation}
\label{invarphiplusgentwo}
{\Xi}^+_{\omega}(f\circ h_t^{(+, \omega)})={\Xi}^+_{\omega}(f).
\end{equation}

We thus have the following


\begin{theorem}
\label{multiplicgen}
Let $\Prob$ be an ergodic $\sigma$-invariant probability measure on $\Omega$ satisfying the 
assumptions \ref{astwo}, \ref{asos}. 
For any $\varepsilon>0$ there exists a positive constant $C_{\varepsilon}$ depending only on 
$\Prob$ such that the following holds. For $\Prob$-almost any $\omega$ there exists a continuous 
mapping $\Xi^+_{\omega}: Lip_w^+(X(\omega))\to {\cal Y}^+_{\omega}$ 
such that for any $f\in Lip_w^+(X(\omega))$, any $x\in X(\omega)$ and all $T>0$
we have
$$
|\int_0^T f\circ h_t^{(+, \omega)}(x) dt - \Xi^+_{\omega}(f)\big(x,t\big)|\leq 
C_{\varepsilon}||f||_{Lip_w^+}(1+T^{\varepsilon}).
$$
The mapping $\Xi^+_{\omega}$ satisfies the equality 
$\Xi^+_{\omega}( f\circ h_t^{(+, \omega)})=\Xi^+_{\omega}(f)$. 
The diagram 
$$
\begin{CD}
Lip_w^+(X(\sigma\omega))@ >\Xi^+_{\sigma\omega}>> {\cal Y}^+_{\sigma\omega} \\
@VVt_{\sigma}^*V           @VVt_{\sigma}^*V   \\
Lip_w^+(X(\omega))@ >\Xi^+_{\omega}>> {\cal Y}^+_{\omega} \\
\end{CD}
$$
is commutative.
\end{theorem}

The mapping $\Xi_{\omega}^+$ is given by $\Xi^+_{\omega}(f)=\Phi_f^+$, where 
$\Phi_f^+$ is defined by (\ref{phifplusgen}).

Now assume that the second Lyapunov exponent $\theta_2$ of the renormalization cocycle ${\mathbb A}$ 
is positive and simple. 
Let $v_2\in E^+_{\omega}$ be a Lyapunov vector 
corresponding to the exponent $\exp(\theta_2)$ (such a vector is 
defined up to multiplication by a scalar).
Introduce a multiplicative cocycle $H^{(2)}(n, \omega)$ over $\sigma$ by the formula 
\begin{equation}
\label{htwo}
H^{(2)}(n, \omega)=\frac{|{\mathbb A}(n, \omega)v_2^{(\omega)}|}{|v_2^{(\omega)}|}.
\end{equation}

Recall that the cocycle $H^{(1)}(n, \omega)$ is given by (\ref{hone}).
Similarly to the above, given a bounded measurable function $f: X(\omega)\to {\mathbb R}$ and $x\in X(\omega)$,
introduce a continuous function ${\mathfrak S}_n[f,x]$ on the unit interval by the formula
\begin{equation}
\label{normsumsgen}
{\mathfrak S}_n[f,x](\tau)=\int\limits_0^{\tau H^{(1)}(n, \omega)} f\circ h^{(+, \omega)}_t(x)dt.
\end{equation}
The functions ${\mathfrak S}_n[f,x]$ are $C[0,1]$-valued random variables on the probability space $(X(\omega), \nu_{\omega})$.

\begin{theorem}
\label{limthmgen}
Let $\Prob$ be an ergodic $\sigma$-invariant probability measure on $\Omega$ satisfying the
assumptions \ref{astwo}, \ref{asos} and such the second Lyapunov exponent of the renormalization 
cocycle ${\mathbb A}$ with respect to $\Prob$ is positive and simple.

For $\Prob$-almost any $\omega^{\prime}\in\Omega$ there exists a non-degenerate compactly supported 
measure $\eta_{\omega^{\prime}}$ on $C[0,1]$ and, for $\Prob$-almost 
any pair $(\omega, \omega^{\prime})$ there exists a sequence of moments 
$l_n=l_n(\omega, \omega^{\prime})$ such that the following holds. 

For $\Prob$-almost any $\omega$ 
there exists a continuous functional 
$$
{\mathfrak a}^{(\omega)}: Lip_w^+(X(\omega))\to {\mathbb R}
$$ 
such that for $\Prob$-almost any $\omega^{\prime}$ and  
any $f\in Lip_{w,0}^+(X(\omega))$ satisfying ${\mathfrak a}^{(\omega)}(f)\neq 0$ 
the sequence of random variables
$$
\frac{{\mathfrak S}_{l_n(\omega,\omega^{\prime})}[f,x]}{{\mathfrak a}^{(\omega)}(f)H^{(2)}(l_n(\omega,\omega^{\prime}), 
\omega)}
$$
converges in distribution to $\eta_{\omega^{\prime}}$ as $n\to\infty$.
\end{theorem}

Theorems \ref{multiplicgen}, \ref{limthmgen} imply Theorems \ref{multiplicmoduli}, \ref{markovlimitmoduli}.
The proofs of Theorems \ref{multiplicgen}, \ref{limthmgen} follow the same pattern 
as those of Theorems \ref{multiplic}, \ref{limthm}; detailed proofs will appear in 
the sequel to this paper.

\end{document}